\documentclass[a4paper,12pt]{article}
\usepackage{mathrsfs}
\usepackage{amssymb}
\usepackage{amsfonts}
\usepackage{latexsym}
\usepackage{amsthm}
\usepackage{amsmath}
\usepackage{pictex}

\def\!{\mskip-\thinmuskip}

\usepackage{graphicx}

\newcommand{\er}[1]{{\rm(\ref{#1})}}
\def\lb{\label}

\theoremstyle{plain}
\newtheorem{theorem}{\bf Theorem}[section]
\newtheorem{lemma}[theorem]{\bf Lemma}

\newtheorem{proposition}[theorem]{\bf Proposition}
\theoremstyle{remark}

\begin{document}

\def\a{\alpha}           \def\cA{\mathcal{A}}     \def\A\mathbf{A}
\def\b{\beta}            \def\cB{\mathcal{B}}     \def\B\mathbf{B}
\def\g{\gamma}           \def\cC{\mathcal{C}}     \def\bC\mathbf{C}
\def\G{\Gamma}           \def\cD{\mathcal{D}}     \def\bD\mathbf{D}
\def\d{\delta}           \def\cE{\mathcal{E}}     \def\bE\mathbf{E}
\def\D{\Delta}           \def\cF{\mathcal{F}}     \def\bF\mathbf{F}
\def\ve{\varepsilon}     \def\cG{\mathcal{G}}     \def\bG\mathbf{G}
\def\z{\zeta}            \def\cH{\mathcal{H}}     \def\bH\mathbf{H}
\def\e{\eta}             \def\cI{\mathcal{I}}     \def\bI\mathbf{I}
\def\vt{\vartheta}       \def\cJ{\mathcal{J}}     \def\bJ\mathbf{J}
\def\vT{\Theta}          \def\cK{\mathcal{K}}     \def\bK\mathbf{K}
\def\k{\kappa}           \def\cL{\mathcal{L}}     \def\bL\mathbf{L}
\def\l{\lambda}          \def\cM{\mathcal{M}}     \def\bM\mathbf{M}
\def\L{\Lambda}          \def\cN{\mathcal{N}}     \def\bN\mathbf{N}
\def\m{\mu}              \def\cO{\mathcal{O}}     \def\bO\mathbf{O}
\def\n{\nu}              \def\cP{\mathcal{P}}     \def\bP\mathbf{P}
\def\r{\rho}             \def\cQ{\mathcal{Q}}     \def\bQ\mathbf{Q}
\def\s{\sigma}           \def\cR{\mathcal{R}}     \def\bR\mathbf{R}
\def\S{\Sigma}           \def\cS{\mathcal{S}}     \def\bS\mathbf{S}
\def\t{\tau}             \def\cT{\mathcal{T}}     \def\bT\mathbf{T}
\def\f{\phi}             \def\cU{\mathcal{U}}     \def\bU\mathbf{U}
\def\F{\Phi}             \def\cV{\mathcal{V}}     \def\bV\mathbf{V}
\def\vp{\varphi}         \def\cW{\mathcal{W}}     \def\bW\mathbf{W}
\def\c{\chi}             \def\cX{\mathcal{X}}     \def\bX\mathbf{X}
\def\p{\psi}             \def\cY{\mathcal{Y}}     \def\bY\mathbf{Y}
\def\P{\Psi}             \def\cZ{\mathcal{Z}}     \def\bZ\mathbf{Z}
\def\o{\omega}
\def\O{\Omega}
\def\x{\xi} \def\X{\Xi}
\def\eps{\epsilon}
\def\vk{\varkappa}

\def\mA{{\mathscr A}}
\def\mB{{\mathscr B}}
\def\mC{{\mathscr C}}
\def\mD{{\mathscr D}}
\def\mE{{\mathscr E}}
\def\mF{{\mathscr F}}
\def\mG{{\mathscr G}}
\def\mH{{\mathscr H}}
\def\mI{{\mathscr I}}
\def\mJ{{\mathscr J}}
\def\mK{{\mathscr K}}
\def\mL{{\mathscr L}}
\def\mM{{\mathscr M}}
\def\mN{{\mathscr N}}
\def\mO{{\mathscr O}}
\def\mP{{\mathscr P}}
\def\mQ{{\mathscr Q}}
\def\mR{{\mathscr R}}
\def\mS{{\mathscr S}}
\def\mT{{\mathscr T}}
\def\mU{{\mathscr U}}
\def\mV{{\mathscr V}}
\def\mW{{\mathscr W}}
\def\mX{{\mathscr X}}
\def\mY{{\mathscr Y}}
\def\mZ{{\mathscr Z}}

\newcommand{\gA}{\mathfrak{A}}
\newcommand{\gB}{\mathfrak{B}}
\newcommand{\gC}{\mathfrak{C}}
\newcommand{\gD}{\mathfrak{D}}
\newcommand{\gE}{\mathfrak{E}}
\newcommand{\gF}{\mathfrak{F}}
\newcommand{\gG}{\mathfrak{G}}
\newcommand{\gH}{\mathfrak{H}}
\newcommand{\gI}{\mathfrak{I}}
\newcommand{\gJ}{\mathfrak{J}}
\newcommand{\gK}{\mathfrak{K}}
\newcommand{\gL}{\mathfrak{L}}
\newcommand{\gM}{\mathfrak{M}}
\newcommand{\gN}{\mathfrak{N}}
\newcommand{\gO}{\mathfrak{O}}
\newcommand{\gP}{\mathfrak{P}}
\newcommand{\gR}{\mathfrak{R}}
\newcommand{\gS}{\mathfrak{S}}
\newcommand{\gT}{\mathfrak{T}}
\newcommand{\gU}{\mathfrak{U}}
\newcommand{\gV}{\mathfrak{V}}
\newcommand{\gW}{\mathfrak{W}}
\newcommand{\gX}{\mathfrak{X}}
\newcommand{\gY}{\mathfrak{Y}}
\newcommand{\gZ}{\mathfrak{Z}}

\def\Z{\mathbb{Z}}
\def\R{\mathbb{R}}
\def\C{\mathbb{C}}
\def\T{\mathbb{T}}
\def\N{\mathbb{N}}
\def\dS{\mathbb{S}}
\def\H{\mathbb{H}}
\def\K{{\Bbb K}}
\def\dD{\mathbb{D}}

\def\qqq{\qquad}
\def\qq{\quad}
\newcommand{\ca}{\begin{cases}}
\newcommand{\ac}{\end{cases}}
\newcommand{\ma}{\begin{pmatrix}}
\newcommand{\am}{\end{pmatrix}}

\def\iint{\int\!\!\!\int}
\def\lt{\biggl}                     \def\rt{\biggr}
\let\ge\geqslant                   \let\le\leqslant
\def\[{\begin{equation}}            \def\]{\end{equation}}
\def\wt{\widetilde}                 \def\pa{\partial}
\def\sm{\setminus}                  \def\es{\emptyset}
\def\no{\noindent}                  \def\ol{\overline}
\def\iy{\infty}                     \def\ev{\equiv}
\def\/{\over}
\def\we{\wedge}
\def\ts{\times}
\def\os{\oplus}
\def\ss{\subset}
\def\h{\hat}
\def\wh{\widehat}
\def\Ra{\Rightarrow}
\def\ra{\rightarrow}
\def\la{\leftarrow}
\def\da{\downarrow}
\def\ua{\uparrow}
\def\lra{\leftrightarrow}
\def\Lra{\Leftrightarrow}
\def\Re{\mathop{\rm Re}\nolimits}
\def\Im{\mathop{\rm Im}\nolimits}
\def\supp{\mathop{\rm supp}\nolimits}
\def\sign{\mathop{\rm sign}\nolimits}
\def\Ran{\mathop{\rm Ran}\nolimits}
\def\Ker{\mathop{\rm Ker}\nolimits}
\def\Tr{\mathop{\rm Tr}\nolimits}
\def\const{\mathop{\rm const}\nolimits}
\def\Wr{\mathop{\rm Wr}\nolimits}
\def\diag{\mathop{\rm diag}\nolimits}
\def\dist{\mathop{\rm dist}\nolimits}

\def\th{\theta}
\def\dlint{\displaystyle\int\limits}
\def\iintt{\mathop{\int\!\!\int\!\!\dots\!\!\int}\limits}
\def\intt{\mathop{\int\int}\limits}
\def\lim{\mathop{\rm lim}\limits}
\def\mult{\!\cdot\!}
\def\BBox{\hspace{1mm}\vrule height6pt width5.5pt depth0pt \hspace{6pt}}
\def\1{1\!\!1}
\newcommand{\bwt}[1]{{\mathop{#1}\limits^{{}_{\,\bf{\sim}}}}\vphantom{#1}}
\newcommand{\bhat}[1]{{\mathop{#1}\limits^{{}_{\,\bf{\wedge}}}}\vphantom{#1}}
\newcommand{\bcheck}[1]{{\mathop{#1}\limits^{{}_{\,\bf{\vee}}}}\vphantom{#1}}
\def\nh{\bhat}
\def\nc{\bcheck}
\newcommand{\oo}[1]{{\mathop{#1}\limits^{\,\circ}}\vphantom{#1}}
\newcommand{\po}[1]{{\mathop{#1}\limits^{\phantom{\circ}}}\vphantom{#1}}
\def\ctg{\mathop{\rm ctg}\nolimits}
\def\notto{\to\!\!\!\!\!\!\!/\,\,\,}

\def\pgbrk{\pagebreak}
\def\Dis{\mathop{\rm Dis}\nolimits}

\def\Twelve{
\font\Tenmsa=msam10 scaled 1200
\font\Sevenmsa=msam7 scaled 1200
\font\Fivemsa=msam5 scaled 1200
\textfont\msbfam=\Tenmsb
\scriptfont\msbfam=\Sevenmsb
\scriptscriptfont\msbfam=\Fivemsb

\font\Teneufm=eufm10 scaled 1200
\font\Seveneufm=eufm7 scaled 1200
\font\Fiveeufm=eufm5 scaled 1200
\textfont\eufmfam=\Teneufm
\scriptfont\eufmfam=\Seveneufm
\scriptscriptfont\eufmfam=\Fiveeufm}

\def\Ten{
\textfont\msafam=\tenmsa
\scriptfont\msafam=\sevenmsa
\scriptscriptfont\msafam=\fivemsa

\textfont\msbfam=\tenmsb
\scriptfont\msbfam=\sevenmsb
\scriptscriptfont\msbfam=\fivemsb

\textfont\eufmfam=\teneufm
\scriptfont\eufmfam=\seveneufm
\scriptscriptfont\eufmfam=\fiveeufm}

\title { Conformal spectral theory for the monodromy matrix
}

\author
{Evgeny Korotyaev
\begin{footnote}
{Institut f\"ur  Mathematik,  Humboldt Universit\"at zu Berlin,
Rudower Chaussee 25, 12489, Berlin, Germany, e-mail:
evgeny@math.hu-berlin.de}
\end{footnote}
}

\maketitle



\begin{abstract}
\no 
For any $N\ts N$ monodromy matrix we define the Lyapunov function,
which is analytic on an associated N-sheeted Riemann surface. On each sheet the Lyapunov function has the standard properties of the Lyapunov
function for the Hill operator. The Lyapunov function has (real or complex) branch points, which we call resonances.
We determine the asymptotics of the periodic, anti-periodic spectrum and of the resonances at high energy. 
We show that the endpoints of each gap are  periodic (anti-periodic) eigenvalues or resonances (real branch points).  Moreover, the
following results are obtained: 1) we define the quasimomentum as
an analytic function on the Riemann surface of the Lyapunov function; various properties and estimates of the quasimomentum are obtained, 2) we construct the conformal mapping with imaginary part given by the
 Lyapunov exponent and we obtain various properties of this
 conformal mapping, which are similar to the case of the Hill operator,
3) we determine various new trace formulae for potentials  and the
 Lyapunov exponent, 4) we obtain a priori estimates of gap
lengths in terms of  the Dirichlet integral.
We apply these results to the Schr\"odinger operators
 and to first order periodic systems  on the real line with a matrix  valued complex self-adjoint periodic potential.

\end{abstract}


\section {Introduction}
\setcounter{equation}{0}

There exist many results about the periodic systems, see 
books \cite{A}, \cite{YS} and interesting papers of Krein \cite{Kr}, Gel'fand and Lidskii \cite{GL},... The basic results for spectral theory for the matrix case were obtained by Lyapunov and Poincar\'e  \cite{YS}. We mention  also new papers \cite{BBK}, \cite{C1}-\cite{C2}, \cite{CK},
\cite{CHGL}, \cite{CG}, \cite{CG1},\cite{GKM},\cite{GS}, \cite{Sh}
and references therin.

Various properties of the Lyapunov function and quasimomentum for scalar 
 Hill operator and the $2\ts 2$ periodic Zakharov-Shabat systems
are well-understood \cite{A}, \cite{M}. Recall the well-known results for the Hill operator $\cS_1 y=-y''+V(t)y$
in $L^2(\R)$ with a periodic potential $V(t+1)=V(t),t\in \R$
and $V\in L^2(0,1)$. The spectrum of $\cS_1$ is purely absolutely continuous and consists of intervals $\wt\s_n=[\l_{n-1}^+,\l_n^-], n\ge 1$. These intervals are separated by the gaps $\g_n=(\l_n^-,\l_n^+)$ of length $|\g_n|\ge 0$. If a gap $\g_n$ is
degenerate, i.e. $|\g_n|=0$, then the corresponding segments
$\wt\s_n,\wt\s_{n+1}$ merge.
The sequence $\l_0^+<\l_1^-\le \l_1^+\ <.....$
is the spectrum of the equation $-y''+Vy=\l y$ with 2-periodic
boundary conditions,  that is  $y(t+2)=y(t), t\in \R$.
 Here equality $\l_n^-= \l_n^+$ means that $\l_n^\pm$ is an eigenvalue of multiplicity 2.  For the equation $-y''+Vy=\l y$ on the real line we define the fundamental solutions $\vt(t,\l)$ and $\vp(t,\l),t\in \R$ satisfying $\vt(0,\l)=\vp'(0,\l)=1, \vt'(0,\l)=\vp(0,\l)=0$. The corresponding monodromy matrix $M$
and the Lyapunov function $\D$  are given by
\[
M(\l)=\ma\vt(1,\l) & \vp(1,\l) \\
\vt'(1,\l) & \vp'(1,\l)\am,\ \ \ \ \
\D(\l)={\Tr M(\l)\/2},\ \ \l\in\C.
\]
 Note that $\D(\l_{n}^{\pm})=(-1)^n, \  n\ge 1$.
The derivative of the Lyapunov function has a zero $\l_n\in [\l^-_n,\l^+_n]$, that is $ \D'(\l_n)=0$ for each $n\ge 1$.
We introduce a conformal mapping (the quasimomentum, 
see \cite{MO}) $k(\cdot):\cZ\to \K$  given by the formula 
$k(z)=\arccos \D (z^2),\ \ z\in \cZ =\C\sm g, \ g=\cup_{n\ne 0} g_n$, where 
$g_n=(z_n^{-},z_n^{+})=-g_{-n}$ and $z_n^{\pm}=\sqrt{\l_n^{\pm }}>0, n\ge 1$ and let $\l_0^+=0$. The quasimomentum domain  
$\K\ev \C \sm\cup_{n\ne 0} c_n$, and vertical slits $c_n=[\pi n+ih_n,\pi n-ih_n]=-c_{-n}$, where a height $h_n\ge 0$  is defined by the equation $\cosh h_n=(-1)^n\D (\l_n)\ge 1$.
 Note that if $V=0$, then $k(z)=z$.
 The following asymptotics hold:
\[
 \lb{T22-1}
k(z)=z-{Q_0\/z}-{Q_2+o(1)\/z^3}\qqq  {\rm  as}\ y>r_0|x|, \qq  y\to\iy,\qq {\rm for \ any} \qq r_0>0,
\]
where $Q_0={1\/2}\int_0^1 V(t)dt={1\/\pi}\int_{g}q(x)dx$ and $Q_2={\|V\|^2\/8}={1\/\pi}\int_{g}x^2q(x)dx$, 
 and $q=\Im k\ge 0$ on $\R$. In particular, this implies
 that if $g=\es$, then $V=0$.

Using the quasimomentum as conformal mapping a priori estimates for various parameters of the Hill operator
and for the Dirac operator  (gap lengths, effective masses,.. in terms of potentials) were obtained in \cite{GT},  \cite{KK1}, \cite{KK2}, \cite{K2}-\cite{K10},\cite{M},\cite{MO}.
Conversly (it is significantly more complicated), a priori estimates of potential in terms spectral data (gap lengths, effective masses,..)
were obtained \cite{K2}, \cite{K5}-\cite{K10}, \cite{M} for example, 
$$
\|V_1\|\le 2\|G\|(1+\|G\|)^{1\/3},\qqq
\|G\|\le 2\|V_1\|(1+\|V_1\|)^{1\/3}, \ \ where \ \ V_1=V-2Q_0,
$$
see [K2], where $\|G\|^2=\sum _{n\ge 1}|\g_n|^2$
and $|\g_n|\ge 0$ is the gap length.
Such a priori  estimates simplify the proof in the inverse spectral
theory for scalar 
 Hill operator \cite{GT}, \cite{KK}, \cite{K1} and for the periodic Zakharov-Shabat systems \cite{K4}, \cite{K5}, see also \cite{K6}, where the author  solved the inverse problem for the operator $-y''+u'y$ on $L^2(\R)$, where periodic $u\in L_{loc}^2(\R)$.

The corresponding theory for the  vector  case is still modest. 
It is well known that the spectrum of
the Schr\"odinger operator on the real line with a $N\ts N$ matrix  valued real periodic potential, $N>1$ is absolutely continuous and consists of intervals separated by gaps \cite{DS}. 
We recall results from \cite{CK} about
this operator:  the Lyapunov function,
which is analytic on an associated N-sheeted Riemann surface is determined. Moreover, the conformal mapping with imaginary part given by the  Lyapunov exponent is constructed and a priori estimates of gap lengths in terms of  potentials are obtained.
Some properties of the monodromy matrices and the corresponding 
the Lyapunov functions for periodic nanotubes  were obtained in
\cite{KL1}, \cite{KL2}.

Introduce the class $\mM_N, N\ge 2$ of monodromy matrices given by

{\bf Definition M.} {\it I) An entire $N\ts N$ matrix-valued function
$M\in \mM_N$ if $M$ satisfies

\no i)  for some unitary matrix $J$  the following identity holds:
\[
\lb{cm-1}
M(z)JM^*(\ol z)=J, \ \ \  \ \  all \ z\in \C.
\]
 $M(z)$ has the eigenvalues $\t_j(z),j\in \ol{1,N}=\{1,..,N\}$ such that:
  
\no ii) If $|\t_j(z)|=1$  for some $(j,z)\in \ol{1,N}\ts\C$,
 then $z\in \R$.

\no II) $M\in \mM_N$ belongs to $\mM_N^0$ 
if each  $\D_j(z)={1\/2}(\t_j(z)+\t_j^{-1}(z)),j\in \ol{N}$ satisfies
\[
\lb{cm-2}
\D_j(z)={1\/2}(\t_j(z)+\t_j^{-1}(z))=\cos z+o(e^{|\Im z|})  \qq as \qq
|z|\to \iy. 
\]
\no III) A matrix-valued function
$M\in \mM_N^0$ belongs to $\mM_N^{0,r},r\ge0$ if
for some constants $C_0,..,C_{r}$ the following asymptotics hold
\[
\lb{ca-3}
 \det (M(z)+M^{-1}(z))= \exp -{iN\rt(z-\sum_0^{2r}{C_s\/z^{s+1}}+{o(1)\/z^{2r+1}}\rt)}\qq  {\rm  as}\ z=iy, \qq  y\to\iy.
\]
}

Note that monodromy matrices for Schr\"odinger operators
(or canonical systems) with periodic matrix- valued potentials
and for Schr\"odinger operators on periodic nanotubes belong to this class, see \cite{A}, \cite{CK}, \cite{KL1}.

The main goal of our paper is to obtain new results
about the Lyapunov functions, the quasimomentum and a priori estimates
of gap lengths in terms of potentials for a class of monodromy
matrices $\mM_N$.
Firstly, we construct Lyapunov functions and  the conformal mapping (averaged quasimomentum) $k(\cdot)$, with imaginary part given by the Lyapunov exponent. 
In fact, we reformulate some spectral problem for the differential operator with periodic matrix coefficients as problems of conformal mapping theory.
Secondly, we obtain various results from the conformal mapping
theory. For solving these "new" problems we use some techniques from [KK1-2], [K2], [K6-8] and [CK]. In particular, we use
the Poisson integral for the domain $\C_+\cup (-1,1)\cup\C_-$. 
We apply these results to the Schr\"odinger operator and to first order periodic systems on the real line with a $N\ts N$ matrix  valued {\bf complex selfadjoint} periodic potential for any $N>1$.
We plan to apply these results to study integrable systems \cite{CD1}, \cite{CD2}, \cite{Ma} and integrated density of states for periodic nanotubes.

An eigenvalue $\t(z)$ of $M(z)$ is called a multiplier.
Note that \er{cm-1} yields that if some $\t(z),z\in\C$ is a multiplier of multiplicity $d\ge 1$, then $1/\ol \t(\ol z)$ 
is a multiplier of multiplicity $d$.

Let $L={M+M^{-1}\/2}$ and $\F(\n,z)=\det (L(z)-\n I_{N})$.
Let $\D_j(z),j\in \ol {1,N}$ be the zeros of $\F(\n,z)=0$,
 where $\ol{m,n}=\{m,m+1,..,n\}$.
 This is an algebraic equation in $\n$ of degree $N$. The coefficients
 of $\F(\n,z)$ are entire in $z\in\C$. It is well known
(see e.g. \cite{Fo}) that the roots $\D_j(z),j\in \ol {1,N}$
constitute one or several branches of one or several analytic functions that have only algebraic singularities in $\C$.
Thus the number of zeros of $\F(\n,z)=0$ is a constant $N_e$
with the exception of some special values of $z$
(see below the definition of a resonance). 
In general, there is an infinite number of such points
on the plane. If all functions $\D_j(z),j\in \ol {1,N}$ are distinct, then $N_e=N$. If some of them are identical,
then $N_e<N$ and $\F(z,\n)=0$ is permanently degenerate.

By definition, the number $z_0$ is a periodic eigenvalue if $z_0$ is a zero of the function $\det (\cM(z)-I_N)$. 
The number $z_1$ is an anti-periodic eigenvalue if
$z_1$ is a zero of the function $\det (\cM(z)+I_N)$. 
We need the following preliminary results

\begin{theorem} \lb{T1}
Let $M\in \mM_N$. Then there exist analytic functions $\wt\D_s, s=1,..,N_0\le N$ on some $N_s$-sheeted Riemann surface $\mR_s, N_s\ge 1$ having the following properties:

\no i) There exist disjoint subsets $\o_s\ss\ol {1,N},
s\in \ol{1,N_0}, \bigcup \o_s=\ol {1,N}$ such that all branches of
$\wt\D_s,s\in \ol{1,N_0}$ are given by $\D_j(z)={1\/2}(\t_j(z)+\t_j^{-1}(z)), \ j\in \o_s$ and satisfy
\[
\lb{T1-1} 
\F(\n,z)=\det \rt({M(z)+M^{-1}(z)\/2}-\n I_{N}\rt)=\prod_1^{N_0} \F_s(\n,z),\qq
 \F_s(\n,z)=\prod_{j\in \o_s}(\n-\D_j(z)),
\]
for any $z,\n\in \C$,
where the functions $\F_s(\n,z)$ are entire with respect to $\n,z\in\C$. Moreover, if $\D_i=\D_j$ for
some $i\in \o_k, j\in \o_s$, then $\F_k=\F_s$ and $\wt \D_k=\wt
\D_s$.

\no ii)  Let 
some branch $\D_j,j\in \ol {1,N}$ be real analytic on some interval
$Y=(\a,\b)\ss\R$ and $-1<\D_j(z)<1$ for any $z\in Y$.
Then $\D_j \ \!\!\!'(z)\ne 0$ for each $z\in Y$ .

\no iii) Each function $\r_s,s=1,..,N_0$, given by \er{T1-2} is entire and  real on the real line,
\[
\lb{T1-2}
\r=\prod_{1}^{N_0}\r_s,\ \ \ 
\r_s(\cdot)=\!\!\!\!\prod_{i<j, i,j\in \o_s}\!\!\!\! (\D_i(\cdot)-\D_j(\cdot))^2.
\]
\no iv) The following identity holds true
\[
\lb{T1-22}
\cup_{j=1}^N \{z\in \C:\D_j(z)\in [-1,1]\}=\R\sm 
\bigcup_{N_-<n<N_+} (z_n^-,z_n^+),\qq  ..<z_n^-<z_n^+
z_{n+1}<..
\]
where $z_n^\pm$ are either periodic (anti-periodic) eigenvalues or real branch points of $\D_j$ (for some $j\in\ol {1,N}$),
which are zeros of $\r$ (below we call such points resonances).

\end{theorem}

{\bf Remark}. 1) If $M\in \mM_N^0$, then  $\r$ is not a polynomial, since $\r$ is bounded on 
$\R$. 


\no 2) Let the surface $\mR=\cup_1^{N_0} \mR_s$ be the union of the disjoint  Riemann surfaces $\mR_s$ and let $\wt\D=\{\wt\D_s, s=1,.,N_0\}$ be the corresponding analytic function on $\mR$.
Let $\z\to z$ be the standard projection from the
surface $\mR$ into the complex plane $\C$. 
{\bf We\  set} $\zeta\in \mR$ and $z=\phi(\zeta)\in \C$. 
The surface $\mR$ is a N-sheeted branch covering of the complex plane,
equipped with the natural projection $\zeta \to z$.
Below we will identify (locally) the point $\zeta\in\mR$ and the point $z=\phi(\zeta)\in \C$ (see \cite{Fo}, 
Chapter 4). In this case we set ${\rm Im}\zeta={\rm Im}\phi(\zeta)$.

\no 3) In the case $M\in \mM_4$ and $\t_1=\t_3, \t_2=\t_4$ 
the function $D(\t,\cdot)=\det (M-\t I_4)$ has the form
 $D(\t,\cdot)=\t^4-T_1\t^3+{1\/2}(T_1^2-T_2)\t^2-T_1\t+1$, and then
\[
\lb{ex2}
 D(\t,\cdot)=\lt(\t^2-2\D_1\t+1\rt)
\lt(\t^2-2\D_2\t+1\rt),\qq 
\D_1={T_1\/2}+{\sqrt{\r}\/2},\qqq \D_2={T_1\/2}-{\sqrt{\r}\/2},
\]
see \cite{BBK}, where $\r={T_2+4\/2}-{T_1^2\/4}$ and $T_j=\Tr M^j(z),\ j=1,2$.

  \vskip 0.25cm
{\bf Definition.} {\it The number $z_0\in \C$ is a {\bf resonance} of $M$, if $z_0$ is a zero of $\r$ given by \er{T1-2}.}

 \vskip 0.25cm

\begin{theorem}   
\lb{Tres} Let $M\in \mM_N$ and let $\D_j, j=\ol{1,\vk}$ have a branch point $z_0\in \R$ for some $\vk\in \N$.
Assume  that each $\D_j(z)\in \R, j\in \ol{1,\vk}$ for all
$z\in (z_0,z_0+\ve)$ (or all for $z\in (z_0-\ve,z_0)$).
Then $\vk=2$ and $z_0$ is a branch point of order ${1\/2}$ for $\D_1, \D_2$
and the function $(\D_1-\D_2)^2$ is analytic in the disk
$\{|z-z_0|<r\}$ for some small $r>0$
and $(\D_1-\D_2)^2$ has a zero $z_0$ of multiplicity $2m+1\ge 1$.
If in addition $\D_1(z_0)\in (-1,1)$, then $m=0$.

\end{theorem}

\no {\bf Remark.} 1) This result is important to describe the spectrum of
Schr\"odinger operators with periodic potentials on armchair
nanotubes \cite{BBKL}. 2) It is very difficult to determine the positions
of resonances. We can use only the Levinson Theorem (see Sect. 2) and Theorem \ref{Tres}. It is similar to another case of poles (other resonances) for scattering  for Schr\"odinger operator with compactly supported potentials
on the real line see \cite{K11}, \cite{Z}.

\no {\bf We consider the conformal mapping associated with  $M\in\mM_N$}. We need functions from the subharmonic  counterpart of the
Cartwright class of the entire functions given by
\[
\lb{dsc} \cS\cC
 =\lt\{q: \C \to \R, \ \begin{array}{c} q {\rm \ is \ subharmonic \ in \ }
 \C {\rm \ and \ harmonic\ outside \ }\R, \cr
 q(\ol z)\equiv q(z),z\in \C,\ \int_{\R}{q_+(t)dt\/1+t^2}<\iy,\
 {\mathop{\lim\sup}\limits}_{z\to\iy}{q(z)\/|z|}<\iy\end{array}\rt\}.
\]
 We recall the class of functions from [KK1] given by
$$
\cS\cK_m^+=\lt\{q\in \cS\cC: q\ge 0, \ \ \lim_{y\to \iy}{q(iy)\/y}=1,\ \ \ 
\int_{\R}(1+t^{2m})q(t)dt<\iy\rt\}, \ \ m\ge 0.
$$
We note that $\cS\cK_{m+1}^+\ss\cS\cK_m^+, m\ge 0$.

Introduce the simple conformal mapping $\e:\C\sm [-1,1]\to \{z\in
\C: |z|>1\}$ by
\[
\lb{de} \e(z)=z+\sqrt{z^2-1},\ \ \ \ z\in \C\sm [-1,1],\ \ {\rm
and}\ \ \ \  \e(z)=2z+o(1)\qq as\ \ |z|\to \iy.
\]
Note that $\e(z)=\ol \e(\ol z), z\in \C\sm [-1,1]$, since $\e(z)>1$
for any $z>1$. The properties of the function $\D_j$ imply $|\e(\wt\D_s(\z))|>1, \z\in \mR_s^+=\{\z\in\mR_s: \Im \z>0\}$. Thus we can define  the quasimomentum $k_j$ (we fix some branch of $\arccos$ and $\D_j(z)$) and the function $q_j$ by
\[
\lb{dkm} k_j(z)=\arccos \D_j(z)=i\log \e(\D_j(z)),\ \  
q_j(z)=\Im k_j(z)=\log |\e(\D_j(z))|,\ \ k\in\ol {1,N} 
\]
 and $z\in \mR_0^+=\C_+\sm \b_+, \b_+=\!\bigcup_{\b\in \cB(\wt\D)\cap \C_+} [\b, \b+i\iy)$,
where $\cB(f)$ is the set of all branch points of the function $f$.
The branch points of $k_j$ in $\C_+$ belong to $\cB(\wt\D)$. 
Define the {\bf averaged quasimomentum} $k$, the {\bf
density} $p$ and the {\bf Lyapunov exponent} $q$ by
\[
\lb{dak} 
k(z)=p(z)+iq(z)={1\/N}\sum_1^N k_j(z), \ \ \  q(z)=\Im k(z),
\ \  z\in \mR_0^+.
\]
Define the sets  
$\s_{(N)}\!=\!\{z\in\R: \D_1(z),..,\D_N(z)\in [-1,1]\}$,
and 
$$
\s_{(1)}\!=\!\{z\in \R: \D_j(z)\in (-1,1),\ \D_s(z)\notin [-1,1] \ {\rm \ some} \ j,s\in\ol {1,N}\}.
$$ 
For the function $k(z)=p(z)+iq(z), z=x+iy\in \ol\C_+$ we formally 
introduce the integrals
\begin{multline}
\lb{dQSI} Q_n={1\/\pi}\int_{\R}x^nq(x)dx,\ \ \ \ \
 I_n^S={1\/\pi}\int_{\R}x^nq(x)dp(x),\ \ \
 I_n^D={1\/\pi}\iint_{\C_+}|\wt k_{(n)}'(z)|^2dxdy, \\
\wt k_{(n)}(z)={1\/\pi}\int_{\R}{t^n q(t)\/t-z}dt =z^n\lt(k(z)-z+\sum_{j=0}^{n-1}Q_sz^{-s-1}\rt), \ \ \
z\in\C_+.
\end{multline}
 Let $C_{us}$ denote the class of all
real upper semi-continuous functions  $h:\R\to \R$. With any $h\in
C_{us}$  we associate the "upper" domain $\K(h)=\{k=p+iq\in\C:
q>h(p), p\in \R\}$. We formulate our first main result. 

\begin{theorem}   
\lb{T2} 
i) Let $M\in \mM_N^0$. Then the averaged quasimomentum $k={1\/N}\sum_1^N k_j$ is analytic
in $\C_+$ and $k:\C_+\to k(\C_+)=\K(h)$ is a conformal mapping
for some $h\in C_{us}$. Moreover, $q=\Im k$ has an harmonic extension
from $\C_+$ into $\O=\C_+\cup \C_-\cup g$ given by $q(z)=q(\ol z), z\in \C_-$ and $q(z)>0$ for any $z\in \O$ and $q\in \cS\cC\cap C(\C)$.
Furhtermore, 
\[
\lb{413} q(z)=y+o(1)\qq {\rm as} \qq |z|\to \iy.\ 
\ \ \
\]
Let in addition, each $\D_j(z)=\cos z+O(|z|^{-{1}}e^{|\Im z|}), j=1,.,N$ as $|z|\to\iy$. Then
\[
\lb{414} q(z)=y+O(|z|^{-{1\/2}}) \qq {\rm as} \qq  |z|\to\iy.
\]

\no ii) Let $M\in \mM_N^{0,r}, r\ge 0$. Then $q\in \cS\cK_{2r}^+$ and there exist branches $k_j,j\in \ol{1,N}$
such that the following asymptotics, identities and estimates hold:
\[
 \lb{T2-1}
k(z)=z-\sum_0^{2r}{Q_s\/z^{s+1}}+{o(1)\/z^{2r+1}}\ \ \  {\rm  as}\ y>r_0|x|, \ \  y\to\iy,\ \
\ \ {\rm for \ any} \ r_0>0,\ 
\]
\[
\lb{T2-2}
C_s=Q_s, \qqq I_s^D+I_{2s}^S=Q_{2s}+{sQ_{s-1}^2\/2}+\sum_{n=0}^{s-2}(n+1)Q_nQ_{2s-2-n},\qq s=0,..,r,
\]
in particular,
\[
 \lb{T2-3}
I_0^D+I_0^S=Q_0, \ \ \  \  if \ \ \  \  r=0 \ \ and \ \ \  \ \ \ I_1^D+I_2^S=Q_2+{Q_0^2\/2} \ \ \  \ 
if \ \ \  \  r=1,
\]
\[
 \lb{T2-4}
q|_{\s_{(N)}}=0,\ \ \ 0<q|_{\s_{(1)}\cup g}\le \sqrt{2Q_0}.\ \ \
\]
\end{theorem}

\no {\bf Remark.} 1) The integral $I_0^S$ is the area between the boundary of $\K(h)$  and the real line.

\begin{figure}[htb]
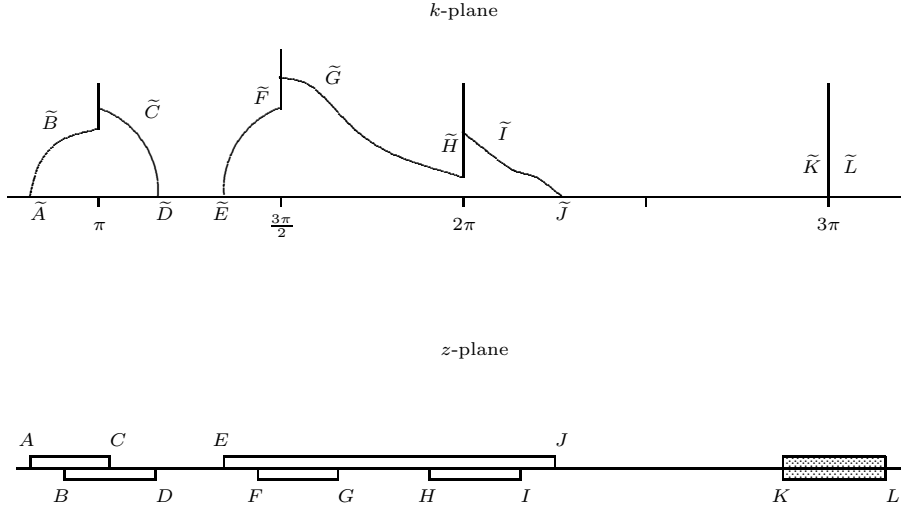

\begin{center}
\font\thinlinefont=cmr5
\begingroup\makeatletter\ifx\SetFigFont\undefined%
\gdef\SetFigFont#1#2#3#4#5{%
  \reset@font\fontsize{#1}{#2pt}%
  \fontfamily{#3}\fontseries{#4}\fontshape{#5}%
  \selectfont}%
\fi\endgroup%
\mbox{\beginpicture \setcoordinatesystem units <0.62992cm,0.62992cm> \unitlength=0.62992cm
\linethickness=1pt \setplotsymbol ({\makebox(0,0)[l]{\tencirc\symbol{'160}}}) \setshadesymbol
({\thinlinefont .}) \setlinear \linethickness= 0.500pt \setplotsymbol ({\thinlinefont .})
{\circulararc 72.861 degrees from  6.638 22.824 center at  7.300 21.129
}%
%
%
\linethickness= 0.500pt \setplotsymbol ({\thinlinefont .}) {\setshadegrid span <1pt>
\shaderectangleson \putrectangle corners at 17.145 15.479 and 19.287 15.240 \setshadegrid span
<5pt> \shaderectanglesoff
}%
%
%
\linethickness= 0.500pt \setplotsymbol ({\thinlinefont .}) {\setshadegrid span <1pt>
\shaderectangleson \putrectangle corners at 17.145 15.240 and 19.287 15.003 \setshadegrid span
<5pt> \shaderectanglesoff
}%
%
%
\linethickness= 0.500pt \setplotsymbol ({\thinlinefont .}) {\putrectangle corners at  5.476
15.479 and 12.383 15.240
}%
%
%
\linethickness= 0.500pt \setplotsymbol ({\thinlinefont .}) {\putrectangle corners at  9.762
15.240 and 11.667 15.003
}%
%
%
\linethickness= 0.500pt \setplotsymbol ({\thinlinefont .}) {\putrectangle corners at  6.191
15.240 and  7.857 15.003
}%
%
%
\linethickness= 0.500pt \setplotsymbol ({\thinlinefont .}) {\putrectangle corners at  2.142
15.240 and  4.047 15.003
}%
%
%
\linethickness= 0.500pt \setplotsymbol ({\thinlinefont .}) {\putrectangle corners at  1.429
15.479 and  3.095 15.240
}%
%
%
\linethickness= 0.500pt \setplotsymbol ({\thinlinefont .}) {\putrule from  1.156 15.240 to
19.6 15.240
}%
%
%
\linethickness= 0.500pt \setplotsymbol ({\thinlinefont .}) {\putrule from  0.953 20.955 to
19.6 20.955
}%
%
%
\linethickness= 0.500pt \setplotsymbol ({\thinlinefont .}) {\putrule from  6.638 23.467 to
6.644 23.467 \plot  6.644 23.467 6.659 23.465 / \plot  6.659 23.465  6.687 23.461 / \plot
6.687 23.461  6.723 23.457 / \plot  6.723 23.457  6.769 23.451 / \plot 6.769 23.451  6.822
23.442 / \plot  6.822 23.442  6.879 23.434 / \plot  6.879 23.434  6.936 23.423 / \plot  6.936
23.423  6.993 23.410 / \plot  6.993 23.410  7.046 23.398 / \plot  7.046 23.398 7.095 23.383 /
\plot  7.095 23.383  7.140 23.368 / \plot  7.140 23.368  7.184 23.351 / \plot  7.184 23.351
7.224 23.332 / \plot 7.224 23.332  7.262 23.309 / \plot  7.262 23.309  7.300 23.285 / \plot
7.300 23.285  7.338 23.258 / \plot  7.338 23.258  7.368 23.235 / \plot  7.368 23.235  7.400
23.207 / \plot  7.400 23.207 7.434 23.180 / \plot  7.434 23.180  7.468 23.148 / \plot  7.468
23.148  7.501 23.116 / \plot  7.501 23.116  7.537 23.080 / \plot 7.537 23.080  7.573 23.044 /
\plot  7.573 23.044  7.612 23.004 / \plot  7.612 23.004  7.650 22.964 / \plot  7.650 22.964
7.688 22.921 / \plot  7.688 22.921  7.728 22.879 / \plot  7.728 22.879 7.766 22.837 / \plot
7.766 22.837  7.806 22.792 / \plot  7.806 22.792  7.846 22.748 / \plot  7.846 22.748  7.885
22.705 / \plot 7.885 22.705  7.925 22.663 / \plot  7.925 22.663  7.963 22.623 / \plot  7.963
22.623  7.999 22.583 / \plot  7.999 22.583  8.037 22.543 / \plot  8.037 22.543  8.073 22.507 /
\plot  8.073 22.507 8.109 22.471 / \plot  8.109 22.471  8.145 22.435 / \plot  8.145 22.435
8.181 22.401 / \plot  8.181 22.401  8.217 22.367 / \plot 8.217 22.367  8.255 22.333 / \plot
8.255 22.333  8.293 22.301 / \plot  8.293 22.301  8.333 22.267 / \plot  8.333 22.267  8.376
22.236 / \plot  8.376 22.236  8.418 22.204 / \plot  8.418 22.204 8.460 22.172 / \plot  8.460
22.172  8.507 22.140 / \plot  8.507 22.140  8.551 22.109 / \plot  8.551 22.109  8.600 22.077 /
\plot 8.600 22.077  8.647 22.047 / \plot  8.647 22.047  8.695 22.020 / \plot  8.695 22.020
8.744 21.990 / \plot  8.744 21.990  8.791 21.965 / \plot  8.791 21.965  8.839 21.937 / \plot
8.839 21.937 8.888 21.912 / \plot  8.888 21.912  8.937 21.888 / \plot  8.937 21.888  8.983
21.865 / \plot  8.983 21.865  9.032 21.844 / \plot 9.032 21.844  9.078 21.823 / \plot  9.078
21.823  9.127 21.802 / \plot  9.127 21.802  9.172 21.785 / \plot  9.172 21.785  9.218 21.766 /
\plot  9.218 21.766  9.267 21.747 / \plot  9.267 21.747 9.318 21.728 / \plot  9.318 21.728
9.370 21.709 / \plot  9.370 21.709  9.426 21.689 / \plot  9.426 21.689  9.487 21.668 / \plot
9.487 21.668  9.553 21.647 / \plot  9.553 21.647  9.622 21.624 / \plot  9.622 21.624  9.696
21.601 / \plot  9.696 21.601  9.777 21.575 / \plot  9.777 21.575  9.859 21.550 / \plot  9.859
21.550 9.946 21.524 / \plot  9.946 21.524 10.031 21.497 / \plot 10.031 21.497 10.116 21.471 /
\plot 10.116 21.471 10.196 21.448 / \plot 10.196 21.448 10.270 21.425 / \plot 10.270 21.425
10.331 21.408 / \plot 10.331 21.408 10.382 21.391 / \plot 10.382 21.391 10.420 21.380 / \plot
10.420 21.380 10.446 21.374 / \plot 10.446 21.374 10.458 21.370 / \plot 10.458 21.370 10.465
21.368 /
}%
%
%
\linethickness= 0.500pt \setplotsymbol ({\thinlinefont .}) {\putrule from  6.668 22.80 to
6.668 24.052
}%
%
%
\linethickness=1pt \setplotsymbol ({\makebox(0,0)[l]{\tencirc\symbol{'160}}}) {\putrule from
18.098 20.955 to 18.098 23.336
}%
%
%
\linethickness= 0.500pt \setplotsymbol ({\thinlinefont .}) {\plot 10.465 22.316 10.471 22.310
/ \plot 10.471 22.310 10.486 22.299 / \plot 10.486 22.299 10.511 22.278 / \plot 10.511 22.278
10.547 22.250 / \plot 10.547 22.250 10.592 22.214 / \plot 10.592 22.214 10.645 22.172 / \plot
10.645 22.172 10.698 22.130 / \plot 10.698 22.130 10.753 22.085 / \plot 10.753 22.085 10.806
22.043 / \plot 10.806 22.043 10.854 22.005 / \plot 10.854 22.005 10.899 21.969 / \plot 10.899
21.969 10.941 21.937 / \plot 10.941 21.937 10.977 21.907 / \plot 10.977 21.907 11.013 21.880 /
\plot 11.013 21.880 11.045 21.855 / \plot 11.045 21.855 11.074 21.831 / \plot 11.074 21.831
11.102 21.808 / \plot 11.102 21.808 11.134 21.783 / \plot 11.134 21.783 11.165 21.759 / \plot
11.165 21.759 11.197 21.734 / \plot 11.197 21.734 11.229 21.711 / \plot 11.229 21.711 11.261
21.685 / \plot 11.261 21.685 11.292 21.662 / \plot 11.292 21.662 11.324 21.639 / \plot 11.324
21.639 11.356 21.618 / \plot 11.356 21.618 11.388 21.596 / \plot 11.388 21.596 11.417 21.577 /
\plot 11.417 21.577 11.447 21.558 / \plot 11.447 21.558 11.474 21.541 / \plot 11.474 21.541
11.502 21.526 / \plot 11.502 21.526 11.527 21.514 / \plot 11.527 21.514 11.553 21.503 / \plot
11.553 21.503 11.576 21.493 / \plot 11.576 21.493 11.604 21.482 / \plot 11.604 21.482 11.633
21.471 / \plot 11.633 21.471 11.663 21.463 / \plot 11.663 21.463 11.692 21.457 / \plot 11.692
21.457 11.722 21.448 / \plot 11.722 21.448 11.756 21.440 / \plot 11.756 21.440 11.788 21.433 /
\plot 11.788 21.433 11.822 21.425 / \plot 11.822 21.425 11.853 21.414 / \plot 11.853 21.414
11.887 21.406 / \plot 11.887 21.406 11.919 21.395 / \plot 11.919 21.395 11.949 21.385 / \plot
11.949 21.385 11.980 21.372 / \plot 11.980 21.372 12.012 21.357 / \plot 12.012 21.357 12.040
21.342 / \plot 12.040 21.342 12.067 21.325 / \plot 12.067 21.325 12.097 21.306 / \plot 12.097
21.306 12.131 21.283 / \plot 12.131 21.283 12.167 21.258 / \plot 12.167 21.258 12.207 21.226 /
\plot 12.207 21.226 12.251 21.192 / \plot 12.251 21.192 12.298 21.154 / \plot 12.298 21.154
12.349 21.114 / \plot 12.349 21.114 12.397 21.074 / \plot 12.397 21.074 12.442 21.035 / \plot
12.442 21.035 12.480 21.004 / \plot 12.480 21.004 12.507 20.983 / \plot 12.507 20.983 12.522
20.968 / \plot 12.522 20.968 12.531 20.961 /
}%
%
%
\linethickness= 0.500pt \setplotsymbol ({\thinlinefont .}) {\putrule from 10.478 21.368 to
10.478 23.336
}%
%
%
\linethickness= 0.500pt \setplotsymbol ({\thinlinefont .}) {\putrule from  2.857 22.384 to
2.855 22.384 \plot  2.855 22.384 2.849 22.382 / \plot  2.849 22.382  2.832 22.377 / \plot
2.832 22.377  2.800 22.371 / \plot  2.800 22.371  2.758 22.360 / \plot 2.758 22.360  2.703
22.348 / \plot  2.703 22.348  2.639 22.333 / \plot  2.639 22.333  2.574 22.316 / \plot  2.574
22.316  2.508 22.299 / \plot  2.508 22.299  2.445 22.280 / \plot  2.445 22.280 2.385 22.265 /
\plot  2.385 22.265  2.330 22.248 / \plot  2.330 22.248  2.282 22.231 / \plot  2.282 22.231
2.237 22.214 / \plot 2.237 22.214  2.197 22.200 / \plot  2.197 22.200  2.159 22.183 / \plot
2.159 22.183  2.125 22.166 / \plot  2.125 22.166  2.091 22.147 / \plot  2.091 22.147  2.057
22.128 / \plot  2.057 22.128 2.026 22.106 / \plot  2.026 22.106  1.994 22.083 / \plot  1.994
22.083  1.962 22.060 / \plot  1.962 22.060  1.930 22.035 / \plot 1.930 22.035  1.901 22.007 /
\plot  1.901 22.007  1.869 21.977 / \plot  1.869 21.977  1.839 21.946 / \plot  1.839 21.946
1.812 21.914 / \plot  1.812 21.914  1.782 21.882 / \plot  1.782 21.882 1.757 21.848 / \plot
1.757 21.848  1.731 21.816 / \plot  1.731 21.816  1.708 21.783 / \plot  1.708 21.783  1.687
21.749 / \plot 1.687 21.749  1.666 21.717 / \plot  1.666 21.717  1.649 21.685 / \plot  1.649
21.685  1.632 21.654 / \plot  1.632 21.654  1.615 21.620 / \plot  1.615 21.620  1.600 21.588 /
\plot  1.600 21.588 1.587 21.556 / \plot  1.587 21.556  1.573 21.520 / \plot  1.573 21.520
1.560 21.484 / \plot  1.560 21.484  1.547 21.444 / \plot 1.547 21.444  1.535 21.402 / \plot
1.535 21.402  1.522 21.355 / \plot  1.522 21.355  1.509 21.304 / \plot  1.509 21.304  1.494
21.249 / \plot  1.494 21.249  1.482 21.194 / \plot  1.482 21.194 1.469 21.137 / \plot  1.469
21.137  1.456 21.086 / \plot  1.456 21.086  1.446 21.040 / \plot  1.446 21.040  1.439 21.002 /
\plot 1.439 21.002  1.433 20.976 / \plot  1.433 20.976  1.431 20.961 / \plot  1.431 20.961
1.429 20.955 /
}%
%
%
\linethickness= 0.500pt \setplotsymbol ({\thinlinefont .}) {\putrule from  2.857 22.384 to
2.857 23.336
}%
%
%
\linethickness= 0.500pt \setplotsymbol ({\thinlinefont .}) {\putrule from 18.098 20.955 to
18.098 20.743
}%
%
%
\linethickness= 0.500pt \setplotsymbol ({\thinlinefont .}) {\putrule from 14.287 20.955 to
14.287 20.743
}%
%
%
\linethickness= 0.500pt \setplotsymbol ({\thinlinefont .}) {\putrule from 10.478 20.955 to
10.478 20.743
}%
%
%
\linethickness= 0.500pt \setplotsymbol ({\thinlinefont .}) {\putrule from  6.668 20.955 to
6.668 20.743
}%
%
%
\linethickness= 0.500pt \setplotsymbol ({\thinlinefont .}) {\putrule from  2.857 20.955 to
2.857 20.743
}%
%
%
\put{\SetFigFont{7}{8.4}{\rmdefault}{\mddefault}{\updefault}{$\wt L$}%
} [lB] at 18.404 21.465
%
%
\put{\SetFigFont{7}{8.4}{\rmdefault}{\mddefault}{\updefault}{$\wt K$}%
} [lB] at 17.554 21.465
%
%
%
%
\put{\SetFigFont{7}{8.4}{\rmdefault}{\mddefault}{\updefault}{$k$-plane}%
} [lB] at  9.762 24.765
%
%
\put{\SetFigFont{7}{8.4}{\rmdefault}{\mddefault}{\updefault}{$z$-plane}%
} [lB] at 10.001 17.621
%
%
%
%
%
%
%
%
%
%
%
%
\put{\SetFigFont{7}{8.4}{\rmdefault}{\mddefault}{\updefault}{$J$}%
} [lB] at 12.383 15.716
%
%
\put{\SetFigFont{7}{8.4}{\rmdefault}{\mddefault}{\updefault}{$I$}%
} [lB] at 11.667 14.527
%
%
\put{\SetFigFont{7}{8.4}{\rmdefault}{\mddefault}{\updefault}{$H$}%
} [lB] at  9.525 14.527
%
%
\put{\SetFigFont{7}{8.4}{\rmdefault}{\mddefault}{\updefault}{$G$}%
} [lB] at  7.857 14.527
%
%
\put{\SetFigFont{7}{8.4}{\rmdefault}{\mddefault}{\updefault}{$F$}%
} [lB] at  5.952 14.527
%
%
\put{\SetFigFont{7}{8.4}{\rmdefault}{\mddefault}{\updefault}{$D$}%
} [lB] at  4.047 14.527
%
%
\put{\SetFigFont{7}{8.4}{\rmdefault}{\mddefault}{\updefault}{$B$}%
} [lB] at  1.905 14.527
%
%
\put{\SetFigFont{7}{8.4}{\rmdefault}{\mddefault}{\updefault}{$A$}%
} [lB] at  1.190 15.716
%
%
\put{\SetFigFont{7}{8.4}{\rmdefault}{\mddefault}{\updefault}{$C$}%
} [lB] at  3.095 15.716
%
%
\put{\SetFigFont{7}{8.4}{\rmdefault}{\mddefault}{\updefault}{$E$}%
} [lB] at  5.239 15.716
%
%
\put{\SetFigFont{7}{8.4}{\rmdefault}{\mddefault}{\updefault}{$L$}%
} [lB] at 19.287 14.527
%
%
\put{\SetFigFont{7}{8.4}{\rmdefault}{\mddefault}{\updefault}{$K$}%
} [lB] at 16.906 14.527
%
%
\put{\SetFigFont{7}{8.4}{\rmdefault}{\mddefault}{\updefault}{$\wt J$}%
} [lB] at 12.383 20.479
%
%
\put{\SetFigFont{7}{8.4}{\rmdefault}{\mddefault}{\updefault}{$\wt I$}%
} [lB] at 11.191 22.147
%
%
\put{\SetFigFont{7}{8.4}{\rmdefault}{\mddefault}{\updefault}{$\wt H$}%
} [lB] at 10.001 21.907
%
%
\put{\SetFigFont{7}{8.4}{\rmdefault}{\mddefault}{\updefault}{$\wt G$}%
} [lB] at  7.586 23.338
%
%
\put{\SetFigFont{7}{8.4}{\rmdefault}{\mddefault}{\updefault}{$\wt F$}%
} [lB] at  6.088 22.9
\put{\SetFigFont{7}{8.4}{\rmdefault}{\mddefault}{\updefault}{$\wt E$}%
} [lB] at  5.239 20.479
%
%
\put{\SetFigFont{7}{8.4}{\rmdefault}{\mddefault}{\updefault}{$\wt C$}%
} [lB] at  3.810 22.623
%
%
\put{\SetFigFont{7}{8.4}{\rmdefault}{\mddefault}{\updefault}{$\wt B$}%
} [lB] at  1.666 22.384
%
%
\put{\SetFigFont{7}{8.4}{\rmdefault}{\mddefault}{\updefault}{$\wt D$}%
} [lB] at  4.047 20.479
%
%
\put{\SetFigFont{7}{8.4}{\rmdefault}{\mddefault}{\updefault}{$\wt A$}%
} [lB] at  1.429 20.479
%
%
\put{\SetFigFont{7}{8.4}{\rmdefault}{\mddefault}{\updefault}{$\pi$}%
} [lB] at  2.75 20.253
%
%
\put{\SetFigFont{7}{8.4}{\rmdefault}{\mddefault}{\updefault}
{$\frac{3\pi}{2}$}%
} [lB] at  6.25 20.253
%
%
\put{\SetFigFont{7}{8.4}{\rmdefault}{\mddefault}{\updefault}{$2\pi$}%
} [lB] at 10.25 20.253
%
%
\put{\SetFigFont{7}{8.4}{\rmdefault}{\mddefault}{\updefault}{$3\pi$}%
} [lB] at 17.858 20.253
%
\linethickness= 0.500pt \setplotsymbol ({\thinlinefont .}) {\circulararc 76.531 degrees from
4.098 20.961 center at  2.290 21.099
}%
\linethickness=0pt \putrectangle corners at  0.906 25.146 and 21.002 14.385
\endpicture}
\end{center}

\caption{$N=2$. The domain $\K(h)=k(\C_+)$ and gaps in the spectrum.
} \label{fig1}
\end{figure}

\no 2) In the case $M\in \mM_4$ and $\t_1=\t_3, \t_2=\t_4$ the Lyapunov functions are  are given by \er{ex2} :  $\D_1={T_1\/2}+{\sqrt{\r}\/2},\  \D_2={T_1\/2}-{\sqrt{\r}\/2}$. 
The mapping $k:\C_+\to \K(h)$ is illustrated in Figure \ref{fig1}.
We have $\D_1^2>1$ on intervals $(A,C), (E,J)$ 
and $\D_2^2>1$ on intervals $(B,D), (F,G), (H,I)$
and $\D_1,\D_2$ are not real on $(K,L)$.
We have also $\wt A=k(A), \wt B=k(B),.., \wt L=k(L)$
and in particular, $k((K,L))=(3\pi, 3\pi +ih_0]$
is a vertical slit for some $h_0>0$.

%
%
%
%
%
%
%
%
%
%
%
%
%
%
%
%
%
%
%
%
%
%
%
%
%

We describe the properties of the conformal mapping $k(\cdot)$.

\begin{theorem}   \lb{T3}
Let $M\in \mM_N^0$. Then the following relations hold:
 \[
 \lb{T3-1}
p_x'(z)\ge 1,\ \ \  z\in\s_{(N)}\ \ {\rm and}\ \ p_x'(z)>0,\ \ \
\ \ z\in\s_{(1)},
\]
here $p_x'(z)=1$ for some $z\in\s_{(N)}$ iff $\s_{(N)}=\s(M)$. Moreover,
\[
\lb{T3-2}
 q_{xx}''(z)<0<q(z),\ \ \ \   p(z)={\rm const}\in {\pi\/N}\Z,
\  {\rm for \ all} \ \ z\in g_n=(z_n^-,z_n^+),
\]
\[
\lb{T3-3} q(x)=q_n^0(x)\rt(1+{1\/\pi}\int_{\R\sm g_n} {q(t)dt\/
q_n^0(t)|t-x|}\rt), \qq
\ \ x\in g_n,\qq q_n^0(z)=|(z-z_n^-)(z_n^+-z)|^{1\/2}, 
\]
\[
 \lb{T3-4}
G^2=\sum_n |g_n|^2\le 8Q_0, \ \ \  and \ \  Q_0\le C_0G^2(1+G^2),
\ \  \ \ if\ \ \ \ \s_{(N)}=\s(M),
\]
for some absolute constant $C_0>0$.
\end{theorem}

 Using this theorem we deduce that the function
$h(p)=q(x(p)), p\in \R$ is continuous on $\R\sm \{p_n, n\in \Z\}$,
where $p_n=p(x), x\in g_n$, and $h(p_n\pm 0)\le h(p_n), \ n  \in \Z$.

{\bf 1. The Schr\"odinger operator}.
We consider the self-adjoint operator $\cS y=-y''+V(t)y,$ acting
in $L^2(\R)^N, N\ge 2$, where $V$ is a 1-periodic $N\ts N$ matrix potential, $V(t)=V^*(t), t\in \R/\Z,$ and $V$ belongs to the complex Hilbert space $\mH$ given by
$$
\mH=\lt\{V(t)=\{V_{jk}(t)\}_{j,k=1}^N,\ \ t\in \R/\Z,\ \ \
\|V\|^2=\int_0^1\Tr V(t)V^*(t)dt<\iy\rt\}.
$$
It is well known (see p.1486-1494 \cite{DS},
 \cite{Ge}) that the spectrum $\s(\cS)$ of $\cS$
is absolutely continuous and consists of non-degenerate
intervals $[\l_{n-1}^+,\l_n^-], n=1,..,N_{G}\le \iy$
and let $\l_0^+=0$.  
These intervals are separated by the gaps
$\g_n=(\l_n^-,\l_n^+)$ with the length $|\g_n|>0$.
 Introduce the fundamental $N\ts N$-matrix solutions $\vp(t,z)$, $\vt(t,z)$ of the equation
\[
\lb{seq}      -f''+Vf=z^2f,\qq \vp(0,z)=\vt'(0,z)=0,
\vp'(0,z)=\vt(0,z)=I_N, \ \ z\in\C,
\]
where $I_N,N\ge 1$ is the identity $N\ts
N$ matrix. Here and below we use the notation
$(')=\pa /\pa t$. We define the monodromy
$2N\ts 2N$-matrix $M$, the matrix $J$ and the trace $T_m, m\in \Z$ by
\[
\lb{dms} M(z)=\ma \vt(1,z)&\vp(1,z)\\
\vt'(1,z)&\vp'(1,z) \am, \ \ J=\ma 0&I_N\\-I_N&0\am ,\ \ \ \ \ \ T_m(z)=\Tr M^m(z).
\]
It is well known  that $M(\cdot)\in\mM_{2N}$,
see \cite{GL}, \cite{Kr}.
 The functions $M$ and $T_m, m\in \Z$ are entire and $\det M=1$. Let $\t_j, j\in \ol{1,2N}$ be the eigenvalues of $M$. 
It is a root of the algebraic equation 
$D(\t,z)\ev\det (M(z)-\t I_{2N})=0, \t,z\in\C$.
Recall that $L={1\/2}(M+M^{-1})$. Each zero of  $\F(\n,z)
=\det (L(z)-\n I_{2N})$ is the Lyapunov function given by 
$
\D_j(z)={1\/2}(\t_j(z)+\t_j^{-1}(z)), \   j\in \ol{1,2N}.
$
For real $V$ we have exactly $\D_j=\D_{N+j},j\in \ol {1,N}$, 
since in this case $\t_{N+j}=1/\t_j$.
Then each $\D_j(z)\in [-1,1]$ gives the spectral point $z^2\in \s(\cS)$ of multiplicity 2, see \cite{CK}. For complex $V$ we have $\D_j,j\in \ol {1,2N}$, where each $\D_j(z)\in [-1,1]$ gives the spectral point $z^2\in \s(\cS)$ of multiplicity 1.
The zeros of $D(1,\sqrt\l)$ ( and $D(-1,\sqrt\l)$) (counted
with multiplicity) are the periodic (anti-periodic) eigenvalues
for the equation $-y''+Vy=\l y$ with periodic (anti-periodic) boundary conditions. 
Let $g=\cup_{n\in \Z} g_n$ where $g_n=(z_n^-,z_n^+), z_n^\pm=\sqrt{\l_n^\pm}>0$
and $g_{-n}=g_n, n\ge 1$. We formulate our main result. 

\begin{theorem}   
\lb{T22}
Let $V=V^*\in \mH$. Then 

\no i) $M$ given by \er{dms} 
belongs to $\mM_{2N}^{0,1}$ and spectrum $\s(\cS)=\R\sm\cup \g_n$, where $\g_0=(-\iy,\l_0^+),\g_n=(\l_n^-,\l_n^+), 1\le n<N_g$ and   
$\l_n^\pm$ are either periodic (anti-periodic) eigenvalues or 
real resonances.

\no ii) The averaged quasimomentum $k={1\/2N}\sum_1^{2N} k_j$ is analytic
in $\C_+$ and $k:\C_+\to k(\C_+)=\K(h)$ is a conformal mapping
onto $\K(h)$ for some $h\in C_{us}$
and $q=\Im k$ has an harmonic extension
from $\C_+$ into $\O=\C_+\cup \C_-\cup g$ given by $q(z)=q(\ol z), z\in \C_-$ and $q(z)>0$ for any $z\in \O$. Furthermore,  $q\in \cS\cK_2^+\cap C(\C)$ and there exist branches $k_j,j\in \ol{1,2N}$
such that the following asymptotics, identities and estimates hold:
\[
 \lb{T22-1}
k(z)=z-{Q_0\/z}-{Q_2+o(1)\/z^3}\qqq  {\rm  as}\ y>r_0|x|, \qq  y\to\iy,\qq {\rm for \ any} \qq r_0>0,
\]
\[
\lb{T22-2} Q_0=I_0^D+I_0^S=\int_0^1 {\Tr V(t)dt\/2N},
\]
\[
\lb{T22-3} Q_2=I_1^D+I_2^S-{Q_0^2\/2}=\int_0^1 {\Tr V^2(t)dt\/8N}=
{\|V\|^2\/8N},
\]
\[
 \lb{T22-4}
q|_{\s_{(N)}}=0,\ \ \ 0<q|_{\s_{(1)}\cup g}\le \sqrt{2Q_0},
\]
 \[
 \lb{T22-5}
k(z)=-\ol{k(-\ol z)},\ z\in\ol\C_+,
\]
\[
 \lb{T22-6}
G^2\ev\sum |\g_n|^2\le {2\|V\|^2\/N},\qqq |\g_n|=\l_n^+-\l_n^-, n\ge 1,
\]
\[
\lb{T22-7}
 \|V\|\le C_0G(1+G^{1\/3}),
\ \  \ \ if\ \ \ \ \s_{(N)}=\s(\cS),
\]
for some absolute constant $C_0$.

\end{theorem}

\no {\bf Remark.} 1) Properties of the Lyapunov functions
are formulated in Theorems \ref{T1}, \ref{Tres}.

\no 2) Various properties of the quasimomnetum $k_j$
are formulated in Theorems \ref{T2}, \ref{T3}, \ref{T41}.

\no 3) The existence of real and complex resonances was proved in [BBK] for the Schr\"odinger operator on the real line with a $2\ts 2$ matrix  real valued periodic potential $V\in \mH$. 

\no 4) If the potential $V\in \mH$ is real and a matrix $\int_0^1V(t)dt$ has distinct eigenvalues,
then the operator  $\cS$ has only finite number complex resonances
[CK].

\no 5) Let $\s(m,A)$ denote the spectrum of a self-adjoint operator  $A$ of multiplicity $m,m\ge 0$.  We have the following simple corollary from Theorem  \ref{T22}: {\it  Let $\s(\cS)=\s(2N,\cS)=\R_+$ for some $V=V^*\in \mH$. Then $V=0$.} There are two simple proofs: 

a) if $\s(\cS)=\s(N,\cS)=\R_+$, then all gaps are close
and the identity \er{T22-3} yields $V=0$.

b) if $\s(\cS)=\s(N,\cS)=\R_+$, then all gaps are close
and the estimate \er{T22-7} yields $V=0$.

\no 6) Recall that the so-called Borg Theorem for periodic systems
was proved in \cite{CHGL},\cite{GKM} for general cases.

{\bf 2. The  periodic  canonical systems}.
Consider the operator $\cK$ given by
$$
\cK y=-iJy'+V(t)y, \qq
  \ J=I_{N_1}\os (-I_{N_2}), \qq  V=\ma 0& v\\
 v^*&0\am,\ \ N_1+N_2=N,\ N_1\ge 1, N_2\ge 1, 
$$
and acting in the space $L^2(\R)^{N}$, where
$v$ is the 1-periodic $N_1\ts N_2$ matrix and $V=V(t)=V(t)^*$ belongs to a subspace $\mH_0\ss\mH$ 
 given by
$$
\mH_0=\lt\{V=\ma \!0\!\!&\! v\!\\ \! v^*&\!0\!\am \in \mH:
v=\{v_{jk}\},\ (j,k)\in \ol{1,N_1}\ts \ol{1,N_2}\rt\}.
$$
It is well known (see [DS] p.1486-1494, [Ge]) that the spectrum
$\s(\cK)$ of $\cK$ is absolutely continuous and consists of
non-degenerated intervals $\s_n$.  These intervals are
separated by the gaps $g_n=(z_n^-,z_n^+)$ with the length
$|g_n|>0, -\iy\le N_g^-< n< N_g^+\le \iy$, where  $N_g=N_g^+- N_g^--1$ is the total number of the gaps.
Introduce the fundamental $N\ts N$-matrix solutions
$\p(t,z)$ of the canonical periodic system
\[
\lb{1e}
-iJ\p'+V(t)\p=z \p,\ \ \ \ z\in\C,\ \ \ \p(0,z)=I_{N}.
\]
It is well known  that the monodromy matrix$\p(1,\cdot)\in \mM_N$, see p.109 \cite{YS}, \cite{Kr}.
The  eigenvalues $\t(z)$  of $\p(1,z)$ are the {\it multiplier} of
$\cK$: to each of them corresponds a solution $f$ of
$\cK \! f=zf$ with $f(t+1)=\t(z)f(t), t\in [0,1)$. They are
roots of the algebraic equation $D(\t,z)\ev\det (\p(1,z)-\t I_{N}), \t,z\in\C$.  The zeros
of $D(1,z)$ (or $D(-1,z)$) are the eigenvalues of  periodic
(anti-periodic) problem for the equation $-iJy'+Vy=zy$.

\begin{theorem}   \lb{T23}
Let $V=V^*\in \mH_0$. Then $\p(1,z)$ belongs to $\mM_{N}^{0,0}$, and
the averaged quasimomentum $k={1\/N}\sum_1^{N} k_j$ is analytic
in $\C_+$ and $k:\C_+\to k(\C_+)=\K(h)$ is a conformal mapping for
some $h\in C_{us}$. Furthermore, and there exist branches $k_j,j\in \ol{1,N}$ such that the following asymptotics, identities and estimates hold:
 \[
 \lb{T23-1}
k(z)-z=-{Q_0+o(1)\/z},\ \ \  if\ y>r|x|, \ \ for \ any \ r>0,
\]
 \[
 \lb{T23-2}
Q_0=I_0^D+I_0^S={\|V\|^2\/2N},
\]
\[
 \lb{T23-4} q|_{\s_{(N)}}=0,\ \ \ 0< q|_{\s_{(1)}\cup g}\le \sqrt{2Q_0},
\]
\[
 \lb{T3-4}
G^2=\sum_n |g_n|^2\le {4\|V\|^2\/N}, 
\]
\[
 \lb{T3-5}
\|V\|\le \sqrt NC_0G(1+G),
\ \  \ \ if\ \ \ \ \s_{(N)}=\s(\cK),
\]
for some absolute constant $C_0$.
\end{theorem}

{\bf Remark.} 1)  In the proof we use arguments from [K4-5]
and [CK]. Theorem \ref{T23} generalizes the result of [CK] to the case of the canonical systems with periodic matrix potential.  

\no 2) We have the following simple corollary from Theorem  \ref{T23}: {\it Let $\s(\cK)=\s(N,\cK)=\R$ for some $V=V^*\in \mH_0$. Then
$V=0$}.  
The  proof is similar to the case of the Schr\"odinger operator,
see Remark 5)  after Theorem  \ref{T22}.

\no 3) In [K3] for
periodic canonical system with  with a specific $4\ts 4$ matrix  real valued symmetric periodic potential the existence of real and complex resonances is proved. 

\no 4) Recall the estimates 
$
{1\/ \sqrt 2}\|G\|\le \|V\|\le 2\|G\|(1+\|G\|)
$
for the case $N_1=N_2=1$ from \cite{K8}.

The plan of our paper is as follows. In Sect. 2 we obtain the basic properties of the Lyapunov functions. In Sect. 3  we obtain the main properties of the quasimomentum and we prove the basic Theorems \ref{T2} and \ref{T3}, devoted to the conformal mapping theory
and Theorems \ref{Tres} about resonances. In Sect. 4 we obtain the results for the Schr\"odinger operator and the first order systems
and Theorem \ref{T22}, \ref{T23} will be proved.

\section {The Lyapunov functions}
\setcounter{equation}{0}

Recall that $D(\t,\cdot)=\det (M-\t I_N)C$ and $L={1\/2}(M+M^{-1})$. We need the simple fact

\begin{proposition}   \lb{T21} i) If $M\in \mM_N$,
then $\F(\n,z)=\det (L(z)-\n I_N),z,\n\in \C$ satisfies
\begin{multline}
\lb{poD} 
 \F(\n,z)=\det (L(z)-\n I_N)=(-1)^N\sum_0^{N}\f_j(z)\n^{N-j},\qq
  \f_0=1,\qq \f_1=-\cT_1,\\ 
\f_2=-{\cT_2+\cT_1\f_1\/2},\qq ..., \qq \f_j=-{1\/j}\sum_1^{j}\cT_k\f_{j-k},..,\qq \f_N=\det L, \qq \cT_m(z)=\Tr L^m(z),
\end{multline}
where each $\f_j$ is  entire and is real on the real line.

\no ii) Let $\D_j(z),j\in \ol {1,N}$ be zeros of the equation
$\F(\n,z)=0$. If some $\D_j(z)\in [-1,1], z\in \C$, then $z\in \R$.

\end{proposition}  
\no{\bf Proof.} i) 
It is well known that the "polynomial" $\F(\n,z), z,\n\in \C$ 
is given by \er{poD}, see p. 331-333 \cite{RS}.  Using the identity \er{cm-1}  we obtain 
$$
\cT_{m}(z)=\Tr L^m(z)={1\/2^m}\Tr \rt(M(z)+M(z)^{-1}\rt)^m
={1\/2^m}\Tr \sum_0^{m}C_p^m M(z)^{2p-m}
$$
$$
={1\/2^m}\sum_0^{m}C_p^m T_{2p-m}(z)
={1\/2^m}\sum_{2p\ge m}^{m}C_p^m \rt(T_{2p-m}(z)+\ol T_{2p-m}(\ol z)\rt), \  \ \ C_m^{N}={N!\/(N-m)!m!}.
$$
This  gives  that each $\cT_m=\Tr L^m, m\in \N$ is the entire
and is real on the real line.

ii) Let some some $\D_j(z)\in [-1,1], z\in \C$. Then
$\t_j(z)$ satifies $|\t_j(z)|=1$ and Definition M yields $z\in \R$.
\BBox

If $\F(\n,z)=(\cos z-\n)^N$, then the Lyapunov function $\D_j(z)=\cos z,
j\in \ol {1,N}$  and the zeros of $\F(1,z)=0$ (or $\F(-1,z)=0$) have the form  $z_{2n,j}=\pi 2n$ (or $z_{2n+1,j}=\pi (2n+1)$ ) for $(n,j)\in \Z\ts \ol {1,N}$. We consider the case $M\in\mM_N^0$.

\begin{lemma} \lb{T31}
Let $M\in \mM_N^0$. Then the following asymptotics hold:
\[
\lb{31} \F(\n,z)-(\cos z-\n)^N=o(e^{N|\Im z|}) \ \ \ \ \ as \ \ \ \ \ |z|\to \iy,
\]
where  $|\n|\le A_0$ for some constant $A_0>0$. Moreover, there exists an integer $n_0$ such that:

\no i) the function $\F(1,z)$ has
exactly $N(2n_0+1)$ roots, counted with multiplicity, in the disc
$\{|z|<\pi(2n_0+1)\}$ and for each $|n|>n_0$, exactly $2N$ roots,
counted with multiplicity, in the domain $\{|z-2\pi
n|<{\pi\/2}\}$. There are no other roots.

\no ii) the function $\F(-1,z)$
has exactly $2Nn_0$ roots, counted with multiplicity, in the disc
$\{|z|<2\pi n_0\}$ and for each $|n|>n_0$, exactly $2N$ roots, counted with multiplicity, in the domain $\{|z-\pi (2n+1)|<{\pi\/2}\}$. There are no other roots.

\no iii)  The function  $D(z,1)$ has only real zeros
$z_{2n,m}, n\in \Z$, their labeling is given by 
\[
\lb{per}
..\le z_{-2,N}\le \underbrace{z_{0,1}\le z_{0,2}\le ...\le z_{0,N}}_{n=0}\le 
\underbrace{z_{2,1}\le ...\le z_{2,N}}_{n=2}\le z_{4,1}\le \dots \ ,
\ \ \ \ \ n\ \ even,
\]
and the function  $D(z,-1)$ has only real zeros
$z_{2n+1,m}, n\in \Z$, their labeling is given by 
\[
\lb{ape} 
..\le z_{-1,N}\le  \underbrace{z_{1,1}\le 
...\le z_{1,N}}_{n=1}\le \underbrace{z_{3,1}\le  ...
\le z_{3,N}}_{n=3}\le z_{5,1}\le \dots �\ , \ \ \ n \ \ odd.
\]
Moreover, they  satisfy
\[
\lb{T1-3}
z_{n,j}=\pi n+o(1)\qqq as\ \ \ \ n\to \pm \iy,\ \ \ \  j\in\ol {1,N} .
\]

\end{lemma} 

\no {\bf Proof.} Each function $\f_n$ in $\F=\sum_0^N \n^{N-n}\f_n$ is the symmetric polynomial of $\D_m, m\in \ol N$. Then $\D_m(z)=\cos z+o(e^{|\Im z|})$ as $|z|\to \iy$ yields \er{31}.

Due to Condition M,ii), the zeros of the function $D(\pm 1,z)$ 
(or $\F(\pm 1,z)$) are real.

  i)   Let $n_1>n_0$ be another integer.
Introduce the contour $C_n(r)=\{z:|z-\pi n|=\pi r\}$. Consider the
contours $C_0(2n_0+1),C_0(2n_1+1),C_{2n}({1\/4}),|n|>n_0$. Then \er{43} and the estimate $e^{{1\/2}|\Im z|}<4|\sin{z\/2}|$ on all contours yield 
\[
\lb{32}
\lt|\F(1,z)-\lt(2\sin{z\/2}\rt)^{2N}\rt|=o(e^{N|\Im z|})=
 o(1) \lt|\sin{z\/2}\rt|^{2N}<{1\/2}\lt|2\sin{z\/2}\rt|^{2N}
\]
for large $n_0$.
Hence, by Rouch\'e's Theorem, $\F(1,z)$ has as many roots, counted
with multiplicities, as $\sin^{2N}{z\/2}$ in each of the bounded
domains and the remaining unbounded domain. Since
$\sin^{2N}{z\/2}$ has exactly one root of the multiplicity $2N$ at
$2\pi n$, and since $n_1>n_0$ can be chosen arbitrarily large, the
point ii) follows.

 ii) The proof for $\F(-1,z)$ is similar.
 
iii) We have $\D_j(z)=\cos z+o(1)$ as   $z=\pi n+O(1)$.
For each $s\in \ol{1,N}$ there exists $j$ such that  $\D_j(z_{n,s})=(-1)^n$. Thus we have $z_{n,s}=\pi n+o(1)$.
$\BBox$

\no {\bf Proof of Theorem  \ref{T1}.} 
The proof repeats the proof of Theorem 1.1 from \cite{CK}.
$\BBox$

We recall some well known facts about entire functions, see [Koo].
An entire function $f(z)$ is said to be of $exponential$ $ type$  if
there is a constant $\g$ such that  $|f(z)|\le $ const. $e^{\g |z|}$
everywhere.  The infimum over the set of $\g$ for which such an inequality holds
is called the type of $f$. The function $f$ is said to belong to
the Cartwright class if $f(z)$ is entire, of exponential type, and
$\int _{\R}{\log ^+|f(x)|dx\/ 1+x^2}<\iy$.
The function $f$ is said to belong to the Cartwright class $\cE_C(a_+,a_-),$
if $f$ is entire, of exponential type, and the following conditions are fulfilled:
$$
\int _{\R}{\log ^+|f(x)|dx\/ 1+x^2}<\iy  ,\ \  \r_\pm(f)=a_\pm,\ \ \
\ \ \ {\rm where}\ \ \
\r_{\pm}(f)\ev \lim \sup_{y\to \iy} {\log |f(\pm iy)|\/y}.
$$

Due to the Paley-Wiener Theorem, if $f\in Cart$ and $f(x)-\cos x\in L^2(\R)$, then $f=\int_{-1}^1\hat \f(t)e^{-itz}dt$ for some $\hat \f\in L^2(-1,1)$. Denote by $\cN (r,f)$ the total number of
zeros of $f$  with modulus $\le r$. Recall the following well
known result (see p.69, [Koo]).

\no    {\bf Theorem (Levinson).}
{\it  Let the function $f\in \cE_C(1,1)$.  Then
$  \cN (r,f)={2\/ \pi }r+o(r)$ as $r\to \iy $.
and for each $\d >0$ the number of zeros of $f$ with modulus $\le r$
lying outside both of the two sectors $|\arg z |, |\arg z -\pi |<\d$
is $o(r)$ for $r\to \iy$.}

We determine asymptotics of
$\r_s, \D_m$ and the zeros of $\r_s, s\in \ol{1,N_0}$.

\begin{lemma} \lb{T32} Let $M\in \mM_N^0$. 
 If $\D_j(z)=\cos z+{b_j\sin z+o(e^{|\Im z|})\/z}$  as $\Im z\to\iy$, and $b_j\neq b_k, j\neq k$ for all $j,k\in \o_s$
for some $s=1,..,N_0$. Then $\r_s\in \cE_C(a,a), a=N_s(N_s-1)$ and
\[
\lb{T32-1}
\r_s(z)=c_s \rt({\sin z\/2z}\rt)^{N_s(N_s-1)}(1+o(1))
\qq {\rm as} \qq \Im z\to\iy, \  \ \ 
c_s=\!\!\prod _{j<k, j,k\in\o_s} \!\!(b_{j}-b_{k})^2.
\]
Let in addition, $\D_j(z)=\cos z+{b_j\sin z+o(e^{|\Im z|})\/z}$  as $ |z|\to\iy$. Then

\no i) There exists an integer $n_0\ge 1$ such that   $\r_s$
has exactly $2N_s(N_s-1)n_0$ roots, counted with multiplicity, in the disc $\{|z|<\pi (n_0+{1\/2})\}$ and for each $|n|>n_0$, exactly $N_s(N_s-1)$ roots, counted
with multiplicity, in the domain $\{|z-\pi n|<{\pi\/2}\}$.
There are no other roots. Moreover, the following asymptotics hold
\[
\lb{T32-2}
\r_s(z)=c_s\rt({\sin z+o(e^{|\Im z|})\/2z}\rt)^{N_s(N_s-1)}
\qq \ {\rm as} \qq  |z|\to\iy.
\]

\no ii) The zeros of $\r_s$ are given by $z_\a^{n\pm}, \a=(j,k), j<k, j,k\in\o_s$ and $n\in \Z\sm \{0\}$. Furthermore, they satisfy
$z_\a^{n\pm}=\pi n+o(1)$ as $n\to \iy$.

\no iii) Let in addition $\D_j(z)=\cos (z-{b_j\/2z})+o(n^{-2})$
 as $|z-\pi n|\le 1, n\to\pm \iy$. Then
\[
\lb{T32-3}
z_\a^{n\pm}=\pi n+{b_j+b_k+o(1)\/2\pi n},\ \ 
\ \ \a=(j,k), \ \ n\to \pm\iy.
\]

\end{lemma} 

\no {\it Proof.} 
  i) using $\D_j(z)-\D_k(z)=(b_{j}-b_{k}){\sin z
\/2z}+o(z^{-1}e^{|\Im z|})$ as $\Im z\to\iy$, we get
$$
\r_s(z)=\prod_{j<k}
(\D_j(z)-\D_k(z))^2=c_s\rt({\sin z+o(e^{|\Im z|})\/2z}\rt)^{N_s(N_s-1)}, \ \ 
$$  
which yields \er{T32-1}.  The proof of \er{T32-2} is similar.

  Let $n_1>n_0$ be another integer.
Introduce the contour $C_n(r)=\{z:|z-\pi n|=\pi r\}$.
 
 Consider the case $N_0=1$, the proof for $N_0\ge 2$ is similar.  
 Let $n_1>n_0$ be another integer. Consider the
contours $C_0(2n_0+1),C_0(2n_1+1),C_{2n}({1\/4}),|n|>n_0$.  Then \er{T32-2} and the estimate
$e^{|\Im z|}<4|\sin z|$ on all contours (for large $n_0$) yield
$
\r(z)=\r^0(z)(1+o(1))$, where $\r^0(z)=c_0({\sin z
\/2z})^{N(N-1)}.
$
Hence, by Rouch\'e's theorem, $\r$ has as many roots, counted
with multiplicities, as $\r^0$ in each of the bounded domains and
the remaining unbounded domain. Since $\r^0$ has exactly one root
of the multiplicity $N(N-1)$ at $\pi n\neq 0$, and since $n_1>n_0$
can be chosen arbitrarily large, the point i) follows.

ii) Thus the zeros of $\r_s$ have the form
$z_\a^{n\pm}, \a=(j,j'), j,j'\in \o_s, j<j', n\in \Z\sm \{0\}$
and satisfy $|z_\a^{n\pm}-\pi n|<\pi/2$. 

We have
$0=\D_j(z)-\D_{j'}(z)=(b_{j}-b_{j'}){\sin z +o(1)
\/2z}$ as $|z_\a^{n\pm}-\pi n|<\pi/2, |n|\to\iy$. Then 
we deduce that $z_\a^{n\pm}=\pi n+o(1)$ as $n\to \iy$.

iii) We have the identity   $\D_j(z)-\D_{j'}(z)=0$ at $z=z_\a^{n\pm}$
 and asymptotics
$$
\cos\lt(z_\a^{n\pm}-{b_j\/2\pi n}\rt)- \cos\lt(z_\a^{n\pm}-{b_{j'}\/2\pi n}\rt)=2(-1)^n\sin {b_{j'}-b_j\/4\pi n}\sin\lt(z_\a^{n\pm}-\pi
n-{b_j+b_{j'}\/4\pi n}\rt)={o(1)\/n^2}
$$
which yields \er{T32-3}.
$\BBox$

\section {The conformal mappings}
\setcounter{equation}{0}

In this section we study  properties of the quasimomentum
and prove theorems about the conformal mapping.
Recall that the Lyapunov function $\wt\D_s(\z)$ is analytic on some 
 $N_s$--sheeted Riemann surface $\mR_s$ and $\mR=\cup_1^{N_0} \mR_s$.
Let $z=x+iy\in\C$ be the natural projection of $\z\in\mR$, $\cB(\wt\D)$ be the set of all branch points of the Lyapunov function and $\mR^\pm=\{\z\in\mR:\pm\Im\z>0\}$. Recall that the simply
connected domains $\mR_0^\pm=\C_\pm\sm\b_\pm\ss\C_\pm$ and define a domain $\mR_0=\C\sm (\b_+\cup \b_-\cup \b_0)$, where
$$
\b_\pm={\bigcup}_{\b\in\cB(\wt\D)\cap\C_\pm}[\b,\b\pm i\infty),
\qqq
\b_0=\{z\in \R: \D_j(z)\notin \R \ for \ some \ j\in \ol {1,N}\}.
$$
Due to the Definition M, $\wt\D(\z)\notin [-1,1]$ for $\z\in\mR^+$. Recall that $q(\z)=|\log\e(\wt\D(\z))|$ is the single-valued on $\mR^+$ imaginary part of the (in general, many-valued on $\mR^+$) quasimomentum $k(\z)=p(\z)+iq(\z)=\arccos\wt\D(\z)=i\log\e(\wt\D(\z))$, where
$$
\e(z)\equiv z+\sqrt{z^2-1},\ \ \ \ \ \e:\C\sm [-1,1]\to \{z\in \C: |z|>1\}.
$$
We denote by $q_j(z)$, $(z,j)\in\C_+\ts \ol {1,N}$, the branches of $q(\z)$ and by $p_j(z)$, $k_j(z)$, $z\in\mR_0^+$, the single-valued branches of $p(\z)$, $k(\z)$, respectively.

\begin{theorem} \lb{T41}
Let $M\in \mM_N$ and let $s\in \ol{1,N_0}$ ($N_0$ is defined in Theorem \ref{T1}). Then
the function $\wt q_s(\z)=\log |\e(\wt\D_s(\z))|$
 is subharmonic on the Riemann surface $\mR_s$. Moreover,

\no 1) If $M\in \mM_N^0$, then the following asymptotics hold
\[
\lb{41} \wt q_s(\z)=y+o(1),\ \ \ \ |\z|\to \iy,\ \ \  \ \Im z\ge 0.
\]
If $\wt\D_s(z)=\cos z+O(1/z)$ as $|z|\to \iy$, then  the following asymptotics hold:
\[
\lb{42} \wt q_s(\z)=y+O(|z|^{-{1\/2}}),\ \ 
|\z|\to \iy, \  \z\in \mR_s, \ \ \Im z\ge 0.
\]


2) Let $\D_j$ be analytic on some bounded interval
$Y=(\a,\b)\ss \R$ for some $j\in \o_s$. Then

\no i) If $\D_j(z)\in \R\sm [-1,1]$ for all $z\in Y$, then $k_j(\cdot)$ has an analytic
extension from $\mR_0^+$ into $\mR_0^+\cup\mR_0^-\cup Y$ such
that
\[
\lb{43} \Re k_j(z)={\rm const}\in \pi \Z,\ \
z\in Y,
\]
\[
\lb{44} q_j(z)=q_j(\ol z)>0,\ \ \ \ z\in \mR_0^+\cup
\mR_0^-\cup Y.
\]
\no ii) If $\D_j(z)\notin \R$ for any $z\in Y$, then there exists a branch $\D_{j'},{j'}\in \o_s$ such that $\ol\D_{j'}(z)=\D_j(z)$ for any
$z\in Y$. The functions $\D_j(z)$ and $k_{j'}+k_j$ have
analytic extensions from $\mR_0^+$ into
$\mR_0^+\cup\mR_0^-\cup Y$ such that
\[
\lb{45}
\D_j(z)=
\ca \D_j(z)\ \ \ &if\ \ \ \ z\in \mR_0^+\\
       \ol\D_{j'}(\ol z)\ \ \ &if\ \ \ \ z\in \mR_0^-\ac,     
\]
\[
\lb{46}
 p_j(z)+p_{j'}(z)={\rm const}\in 2\pi \Z,\ \
z\in Y,\ \ \
q_j(z)=q_{j'}(z),\ \ \ \ \ z\in Y,
\]
\[
\lb{47}
 q_j(z)+q_{j'}(z)=q_j(\ol z)+q_{j'}(\ol z)>0,\ \
\ z\in \mR_0^+\cup \mR_0^-\cup Y.
\]
\end{theorem}
\no {\bf Proof.} 
The proof repeats the proof of Theorem 4.1 from \cite{CK}.
$\BBox$

{\bf Proof of Theorem \ref{Tres}}. 
Puiseux series for $\D_j,j=1,..,\vk $ are given by
\[
\lb{br}
\D_j(z)=\D_1(z_0)+a_1e t+a_2e^2t^2+..,\qq e=e^{i{2\pi\/p}},\qq  t=(z-z_0)^{1\/\vk}\in D_r=\{t:|t|<r\}.
\]
If $j=\vk$, then we obtain $a_n\in \R$ for all $n\ge 1$.
Furthermore, $\vk=2$, since $\D_j(z)$ is real for all $j\in \ol{1,\vk}$ and all $z\in (z_0,z_0+\ve)$ (or all $z\in (z_0-\ve,z_0)$).
Moreover, these arguments give
\[
\lb{br1}
\D_j(z_0+t^2)=f(t)+(-1)^jt^{1+2m}g(t), \qq j=1,2,\qq t\in D_r, \ m-1\in \N, \ g(0)\ne 0,
\]
where $f,g$ are analytic function in the disk $D_r$ for some 
$r>0$ and $f,g$ are real on $(-r,r)$.
Thus $F=(\D_1-\D_2)^2$ is analytic in $D_r$ and
$F(z)=z^{1+2m}(C+O(z))$ as $z\to z_0$ for some $C\ne 0$.

Let $\D_1(z_0)\in (-1,1)$. Then using \er{br1}, identity $
\D_j(z)=\cos k_j(z)$ and the Implicit Function Theorem, we obtain
$$
k_j(z_0+t^2)=f_j(t)+t^{1+2m_j}g_j(t), \qq j=1,2,\qq t\in D_r, \ m_j-1\in \N, \ g_j(0)\ne 0,
$$
where $f_j,g_j$ are some analytic functions in the disk $D_r$ for some $r>0$. Then \er{47}  yields
$$
k_1(z_0+t^2)+k_2(z_0+t^2)=2\pi n_0+t^s(c+O(t))\qq as \qq t\to 0,
t\in\ol\C_+,
$$
for some $(n_0,s)\in \Z\ts\N$ and $c\ne 0$. The case $s\ge 2$ is absent,
since $q_1,q_1>0$ in $D_r\cup \C_+$ for some $r>0$, which yields $s=m_1=m_2=1$.
\BBox

Recall the needed properties of the functions $q\in\cS\cC$
defined in Sect. 1 and $k=p+iq$. It is well known, that $p\in C(\ol\C_+)$ and
 ${1\/2\pi}(\pa_x^2+\pa_y^2) q=\m_q$ (in a sense of
distribution) is a so-called Riesz measure of the function $v$.
Moreover, the following identities are fulfilled:
\[
\lb{48} \pi\m_q((x_1,x_2))=p(x_2)-p(x_1), \ \ \ \ {\rm for \ any
}\ \ x_1<x_2,\  x_1,x_2\in \R,
\]
\[
\lb{49} {\pa q(z)\/\pa y}=y\int_{\R}{d\m_v(t)\/(t-x)^2+y^2},\ \ \
z=x+iy\in\C_+,
\]
which yields ${\pa q(z)\/\pa y}\ge 0, z\in\C_+$. Moreover,
$q(x)=q(x\pm i0), x\in \R.$ It is well known that if $q\in
\cS\cC$, then
\[
\lb{410} \int_{\R}{d\m_q(t)\/1+t^2}<+\iy,\ \ \
\lim\sup_{z\to\iy}{q(z)\/|z|}=
\lim_{y\to+\iy}{k(iy)\/iy}=\lim_{y\to+\iy}{q(iy)\/y}=\lim_{x\to
\pm \iy}{p(x)\/x}\ge 0.
\]
 Now we recall the well known fact (see [Ah]).

\no {\bf Theorem (Nevanlinna).}
 {\it \no i) Let $\m$ be a Borel measure on
$\R$ such that  $\int_{\R}(1+x^{2r})d\m(x)<\iy$ for some integer
$r\ge 0$. Then for each $s>0$ the following asymptotics hold
$$
\int_{\R}{d\m(t)\/t-z}=-\sum_{k=0}^{2r}{Q_k\/z^{k+1}}+
{o(1)\/z^{2p+1}},\ \ \ |z|\to\iy,\ y>s|x|,\ \ \ \
Q_n=\int_{\R}x^nd\m(x),\ 0\le n\le 2r.
$$
\no ii) Let $f$ be an analytic function in $\C_+$ such that $\Im
f(z)\ge 0$ for all $z\in\C_+$ and
\[
\lb{411}
\Im f(iy)=c_0y^{-1}\!+\!\dots\!+\!c_{2p-1}y^{-2r}\!+\!O(y^{-2r-1})
\ \ \ {\rm as}\ \ y\!\to\!\iy
\]
for some $c_0,..c_{2r-1}\in \R$ and $r\ge 0$.  Then $
f(z)=C+\int_{\R}{d\m(t)\/t-z},\ z\in\C_+$, for some Borel measure
$\m$ on $\R$ such that $\int_{\R}(1+x^{2r})d\m(x)<\iy$ and $C\in
\R$.}

\no {\bf Proof of Theorem \ref{T2}.} 
i) The proof repeats the proof of Theorem 4.2 from \cite{CK}.

ii) Below for each real harmonic  function $q(z),z\in\C_+$ we introduce an analytic function $k=p+iq$  in $\C_+$, where $(-p)$ is some harmonic conjugate of $q$ for $\C_+$. If
$q\in \cS\cK_0^+$, then the function $k=p+iq$ in $\C_+$ is defined
by
\[
\lb{51} k(z)=z+{1\/\pi}\int_{\R}{q(t)dt\/t-z},\ \ \ \ z\in\C_+.
\]

Due to i), the
function $q={1\/N}\sum_1^{N}q_j\in \cS\cC\cap C(\C)$ and $q$ is positive
in $\C_+$. Let  $z=iy, y\to \iy$, using   $q_m(iy)=y+o(1)$ (see \er{41} ) we obtain
\[
\F(z,0)=\prod_1^N
\D_m(z)={(1+e^{-2Ny}O(1))\/2^N}e^{-i\sum_1^N k_m(z)}
={(1+e^{-2Ny}O(1))\/2^N}e^{-iNk(z)}.
\]
Thus these asymptotics and  
$
\F(z,0)=(\cos^{N}
z)\exp i\rt({-\sum_{j=0}^{2r}{C_j\/z^{j+1}}+{o(1)\/z^{2r+1}}}\rt)
$, see \er{ca-3}, give
\[
k(z)=z-\sum_{j=0}^{2r}{C_j\/z^{j+1}}+{o(1)\/z^{2r+1}}\qqq \ as  \qq 
z=iy,\ y\to \iy.
\]
 Then by i) and  the Nevanlinna Theorem,
 the function $q\in \cS\cK_{2r}^+$.  
 We need the following  result from [KK1]:
 {\it Let $q\in \cS\cK_m^+$ for some $m\ge 0$ and $q\neq const$. Then $k:\C_+\to k(\C_+)=\K(h)$
is a conformal mapping for some $h\in C_{us}, h\ge 0$. Moreover, the
following asymptotics, estimates and identities are fulfilled:
\[
\lb{52} 
k(z)=z-\sum_{k=0}^{2m}{Q_k\/z^{k+1}}+{o(1)\/z^{2m+1}}\ \ as \ \ |z|\to \iy,\ \ y\ge r|x|,
\ \ {\rm for \ any}\ r>0,
\]
\[
\lb{53}  I_s^D+I_{2s}^S=Q_{2s}+{sQ_{s-1}^2\/2}+\sum_{n=0}^{s-2}(n+1)Q_nQ_{2s-2-n}, \ \ \ \ \ \ s=0,...,m,
\]
\[
\lb{54} 
\sup_{x\in\R}q^2(x)\le 2Q_0,
\]
where $I_n^S, I_n^D, Q_n, n\ge 0$ are given by \er{dQSI}.}
 Thus the above results give that  $k:\C_+\to k(\C_+)=\K(h)$ is a conformal mapping for some $h\in C_{us}$. Moreover, 
 asymptotics \er{T2-1} and identities \er{T2-2}, \er{T2-3}
hold true.

Estimate \er{54} gives $q |_\R\le\sqrt{2Q_0}$.
If $z\in \s_{(N)}$, then $q(z)=0$, since $\D_j(z)\in[-1,1]$
 and $q_j(z)=0$ for all $j=1,..,N$. 
If $z\in \s_{(1)}\cup g$, then  $\D_j(z)\notin [-1,1]$ for some $j=1,..,N$. Thus $q_j(z)>0$ and $q(z)>0$,  
which gives \er{T2-4}. \BBox

\no {\bf Proof of Theorem \ref{T3}.} i) We show \er{T3-1}. Using
$k(z)=z+{1\/\pi}\int_{\R}{q(t)dt\/(t-x)}, z\in \C_+$ and  \er{T2-4}  we obtain
$k'(z)=1+{1\/\pi}\int_{\R}{q(t)dt\/(t-x)^2}>1,\ \ z\in \s_{(N)}$,
which yields $p_x'(z)>1,z\in \s_{(N)}$.
Moreover, we have $p_x'(z)=1$ for some $z\in \s_{(N)}$ iff $q|_\R=0$.

The function $\wt q_s$ is harmonic in $\mR_s^+=\{\z\in\mR_s: \Im \z>0\}$. The function $f=\wt q_s-y$ is harmonic in $\mR_s^+$
and $f>0$.

Let $x\in \s_{(1)}$.  Then some branch $\D_m(x)\in (-1,1)$ for all $x\in Y\ss\s_{(1)}$ for some small interval  $Y=(\a,\b)$.  We have $\D_m(x)=\cos
k_m(x)$. Thus we get $k_m'(x)=-{\D_m'(x)\/\sin
k_m(x)}\ne 0$. 
Hence we get $p_m'(x)=k_m'(x)>0$, since $q_m(x)=0$ and \er{49}
yields $p_x(x)\ge 0$ for $x\in Y$, which gives \er{T3-1}.

We show \er{T3-2}. 
 Note that $q(z)>0,z\in g$, otherwise we have not a gap. Using \er{49} we obtain $ -{\pa^2 q(z)\/\pa x^2}
={\pa^2 q(z)\/\pa y^2}=\int_{\R}{d\m_q(t)\/(t-x)^2}>0,\ \ \ z=x\in
g_n, $ which yields $q_{xx}''(z)<0, z\in g_n$. Consider the function $p(z), z\in g_n$.
Theorem \ref{T41} yields
$
N\Re k(z)=\sum_1^N \Re k_m=\sum_1^N \pi n_m=\pi N_{n}, 
$
which gives \er{T3-2}

We recall the result from \cite{KK2}. Let a function $f$ be harmonic
and positive in the domain $\C\sm \ol g_n, g_n=(z_n^-,z_n^+)\ne \es$ and
$f(iy)=y(1+o(1))$ as $y\to\iy$. Assume $f(z)=f(\ol z),\ z\in
\C\sm \ol g_n$ and $f\in C(\ol{\C_+})$. Then
$$
f(x)=q_n(x)\lt(1+{1\/\pi}\int_{\R\sm I} {f(t)dt\/
|t-x|q_n(t)}\rt),\ \  \ q_n(x)=|(z-z_n^-)(z_n^+-z)|^{1\/2},\qq
 x\in g_n.
$$
Hence the last identity and properties of $q$ yield \er{T3-3}
and the estimate $q(z)\ge q_n(z),\ z\in g_n=(z_n^-,z_n^+)$. 
This estimate implies $Q_0\ge {1\/\pi}\sum\int_{g_n} q_n(t)dt= {1\/8}\sum |g_n|^2$, which yields the first estimate in \er{T3-4}.
\BBox

\section {Proof of Theorems \ref{T22} and \ref{T23}}
\setcounter{equation}{0}

We begin with some notational
convention. A vector $h=\{h_n\}_1^N\in \C^N$ has the Euclidean
norm $|h|^2=\sum_1^N|h_n|^2$, while a $N\ts N$ matrix $A$ has the
operator norm given by $|A|=\sup_{|h|=1} |Ah|$. Note that
$|A|^2\le \Tr A^*A$. 

{\bf 1  The first order periodic systems.}  In this case $J=I_{N_1}\os (-I_{N_2})$. Below we use arguments from \cite{K4},\cite{K5}. We need the identities 
\[
\lb{61} JV=-VJ,  \ \ \  e^{zJ} V=Ve^{-zJ},\ \ \ all\qq
(z,V)\in \C\ts\mH_0.  
\]
 The solution of the equation
$-iJ\p'+V\p=z\p, \p_0(z)=I_N$ satisfies the integral equation
\[
\lb{62} \p(t,z)=e^{iztJ}-i\int_0^te^{izJ(t-s)}JV(s)\p(s,z)ds,\ \
\ t\ge 0,\ z\in \C,
\]
and $\p$ is given by
\[
\lb{63} \p(t,z)=\sum_{n\ge0}\p_n(t,z), \ \ \
\p_n(t,z)=-i\int_0^te^{izJ (t-s)}J V(s)\p_{n-1}(t,z)ds,\ \
\ \p_0(t,z)=e^{iztJ},
\]
$\  n\ge1$.  Using \er{61}, \er{63} we have
\[
\lb{64} \p_1(t,z)=-i\int_0^te^{izJ (t-s)}J V(s)e^{izsJ }ds=
-i\int_0^te^{izJ (t-2s)}J V(s)ds,
\]
\[
\lb{65} \p_2(t,z)=-i\int_0^te^{izJ (t-t_1)}J V(t_1)\p_1(t_1,z)dt_1=
\int_0^tdt_1\int _0^{t_1}e^{izJ(t-2t_1+2t_2)}V(t_1)V(t_2)dt_2.
\]
Proceeding by induction,
\[
\lb{66} \p_{2n}(1,z)=\!\!\int_0^1\!\!\!\!dt_1\!\!\dots \!\!\int _0^{t_{2n-1}}\!\!\!\!
e^{izJ\vk_{2n}}V(t_1)\dots V(t_{2n})dt_{2n},\ \ \
\vk_n=1-2t_1+2t_2\dots +(-1)^n2t_n,
\]
\[
\lb{67} \p_{2n+1}(1,z)=-iJ\int_0^1dt_1\dots \int _0^{t_{2n}}
e^{izJ\vk_{2n+1}}V(t_1)\dots V(t_{2n+1})dt_{2n+1}.
\]
 We need the following estimates.

\no  \begin{lemma}\lb{T61} 
Let $V\in\mH_0$.  For each $z\in \C$ there
exists a unique solution $\p$ of Eq. \er{62} given by \er{63} and
series \er{63} converge uniformly on bounded subsets of
$\R\ts\C\ts \mH$. For each $t\ge 0$ the function $\p (t,z)$ is
entire on $\C$.  Moreover, for any $n\ge 0$ and $(t,z)\in
[0,\iy)\ts \C$ the following estimates and asymptotics hold:
\[
\lb{68} |\p_n(t,z)|\le {e^{|\Im z|t}\/n!}\lt(\int_0^t |V(s)|ds
\rt)^n,
\]
\[
\lb{69} |\p(t,z)-\sum_0^{n-1}\p_j(t,z)|\le {(\sqrt
t\|V\|)^n\/n!}e^{t|\Im z|+\int_0^t |V(s)|ds},
\]
\[
\lb{610} \p(t,z)-e^{iztJ}=o(e^{t|\Im z|})\ \ \ \ as  \ \
\ |z|\to \iy,
\]
\[
\lb{611}
T_m(z)=N_1e^{iz}+N_2e^{-iz} +o(e^{m|\Im z|})\ \ as  \ \ \ |z|\to \iy.
\] 
If $V^{\n},V\in \mH_0$ and if the sequence $V^{\n}\to V$ weakly in $\mH_0$ as $\n\to \iy$, then $\p (t,z,V^{\n})\to \p (t,z,V)$
 uniformly on bounded subsets of $\R\ts \C$.
\end{lemma}
\no {\it  Proof.}  \er{66}, \er{67} give
$$
 |\p_n(t,z)|\le \!\!\int_0^t\!\!dt_1\int _0^{t_1}\!\!dt_2\dots \int
_0^{t_{n-1}}\!\! e^{|\Im z|t}|V(t_1)|\dots |V(t_n)|dt_n\le 
{e^{|\Im z|t}\/n!}\lt(\int_0^t |V(s)|ds
\rt)^n,
$$
since $|t-2t_1+2t_2\dots +(-1)^n2t_n|\le t$, which yields \er{68}. This shows that for each $t
>0$ the series \er{63} converges uniformly on bounded subsets of
$\C\ts\cH_0$. Each term of this series is an entire function.
Hence the sum is an entire function. Summing the majorants we
obtain estimates \er{69}. The proof of asymptotics
in \er{610} is standard (see e.g. \cite{K4} or \cite{K5}).
\er{610} implies \er{611}.

Assume that the sequence $V^{\n}\to V$ weakly in $\mH$, as $\n\to \iy$. Then each term $\p_n(t,z,V^{\n})\to \p_n(t,z,V)$
 uniformly on bounded subsets of $\R\ts \C$ and fixed $n\ge 1$.
Then \er{69} gives that $\p (t,z,V^{\n})\to \p (t,z,V)$
 uniformly on bounded subsets of $\R\ts \C$.
 $\BBox$

We need asymptotics of $L(z)={1\/2}(M(z)+M^{-1}(z))$ as $\Im z\to \iy$. 
The identities \er{cm-1}, \er{63} yield
\[
\lb{612} L(z)={1\/2}(\p(1,z)+J\p^*(1,\ol z)J)=\sum_{n\ge 0} L_n(z),\
\ L_n(z)={1\/2}(\p_n(1,z)+J\p_n^*(1,\ol z)J),
\]
where series \er{612} converge uniformly on bounded subsets of
$\C\ts \mH_0$. Using \er{61},\er{66}, \er{67}, we obtain
\[
\lb{613} J{\p_{2n}}^*(1,\ol z)J=\int_0^1dt_1\dots \int _0^{t_{2n-1}}
e^{-izJ\vk_{2n}}V({t_{2n}})\dots V({t_1})dt_{2n},
\ \ \ \  \ 
\]
\[
\lb{614} J{\p_{2n+1}}^*(1,\ol z)J=iJ\int_0^1dt_1\dots \int _0^{t_{2n}}
e^{izJ\vk_{2n+1}} V({t_{2n+1}})\dots
V({t_1})dt_{2n+1},
\]
and in particular, $J\p_1^*(1,\ol z)J=-\p_1(1,z)$.
Thus we deduce that
\[
\lb{615} L=\cos z+L_2+\sum_{n\ge 3} L_n,\ \ L_1=0,
\]
\[
\lb{616}
  L_2(z)={1\/2}\int_0^1\!\!\!\int
_0^{t_1}\rt(e^{izJ(1-2t_1+2t_2)}V({t_1})V({t_2})
+e^{-izJ(1-2t_1+2t_2)}V({t_2})V({t_1})\rt)dt,....
\]

\begin{lemma}  \lb{T62} 
Let $V\in \mH_0$. Then for $y\ge r_0|x|, y\to \iy,r_0>0 $ the following asymptotics hold: 
\[
\lb{617} 
\int_0^1\!\!dt_1\!\!\dots \!\!\int _0^{t_{n-1}}\!\!\!\! e^{\pm izJ\vk_{n}}V({t_1})\dots V({t_{n}})dt_{n}
={o(e^{y})\/|z|},
\]
\[
\lb{618} L_2(z)={\sin z\/2z }(\mV+o(1)),\ \ \ \ where\ \ \ 
\mV=\int_0^1V^2(t)dt,
\]
\[
\lb{619} L(z)=\cos \rt(zI_{N}-{\mV+o(1)\/2z }\rt),
\]
\[
\lb{620}
\det L(z)=(\cos^{N} z)\exp \rt({i\|V\|^2+o(1)\/2z}\rt).
\]
\end{lemma}

 \no {\it Proof.} Let $V_t=V(t)$. Asymptotics \er{TA-1}\er{TA-2} give
\[
\lb{621} \int_0^1dt\int_0^te^{i2z(t-s)}V_tV_sds
={i\/2z}\rt(\mV+o(1)\rt),\ \ \ 
 \int_0^1dt\int_0^t
e^{-i2z(t-s)}V_tV_sds={o(e^{2y})\/|z|}.
\]
 Using \er{621} and $\vk_2=1-2t+2s$ we get
 \[
\lb{622} L_2(z)={1\/2}\int_0^1\!\!\!\int
_0^{t}\rt(e^{izJ\vk_2}V_{t}V_{s}
+e^{-izJ\vk_2}V_{s}V_{t}\rt)dt
\]
$$
={1\/2}\int_0^1\!\!\!\int
_0^{t_1}\rt(\cos z\vk_2(V_{t}V_{s}+V_{s}V_{t})
+iJ\sin z\vk_2 (V_{t}V_{s}-V_{s}V_{t})\rt)dt={i\cos\/2z}\rt(\int_0^1 V_t^2dt+o(1)\rt),
$$
which gives \er{618}.  The similar arguments yield \er{617}. 
 
In order to estimate $L_n$ we consider the function $f_n^\pm=f_n^\pm(y)=\int_{G_n}\!\! e^{\pm y\vk_n}dt$. 
Recall $\vk_n=1-2t_1+2t_2\dots +(-1)^n2t_n,\ t=(t_1,..,t_n)\in  
G_n=\{0< t_n<... t_2<t_1<1\}\ss \R^n$. We need the simple estimate
\[
\lb{623}
\int_0^a e^{2yt}dt\le {e^{2ya}\/2y},\  \ \  \ \ \ 
\int_0^a\!\!\!\int_0^{t}e^{-2y(t-s)}dtds\le {1\/2y}, \qq \ 
any \qq a,y>0.
\]
The direct calculations imply
$
f_p^\pm \le e^yy^{-2},y\ge 1,\  p=2,3.
$
Consider $f_{n}^+$, the proof for $f_{n}^-$ is similar.  Using the second estimate from \er{623} we have
$
f_{2n}^+(y)\le  {e^{y}\/(2y)^n}
$
This and the first estimate from \er{623}  yield
$
f_{2n+1}^+(y)\le  {f_{2n}^+(y)\/2y}\le {e^{y}\/(2y)^{n+1}}.
$
Thus we obtain
\[
\lb{624}
f_{n}^\pm(y)\le  {e^{y}\/(2y)^{n\/2}},\ \ \ n,y\ge 1.
\]
Hence we have
\[
\lb{625}
|L_n(z)|\le \!\!\int_{G_n}\!\! (e^{y\vk_n}+e^{-y\vk_n})\prod_1^n
|V(t_j)| dt \le \rt(\int_{G_n}\!\! (e^{y\vk_n}+e^{-y\vk_n})^2
dt\rt)^{1\/2} \ {\|V\|^{n}\/\sqrt{n!}}\le {2e^{y}\/\sqrt{n!}}\ \ve^n,
\]
where $\ve ={\|V\|\/(4y)^{1\/4}}\to 0$, since 
$$
\int_{G_n}\!\! \prod_1^n |V_{t_j}|^2 dt={\|V\|^{2n}\/n!},\ \ \ \ 
\!\!\int_{G_n}\!\! (e^{y\vk_n}+e^{-y\vk_n})^2dt\le
2\!\!\int_{G_n}\!\! (e^{2y\vk_n}+e^{-2y\vk_n})^2dt\le {4e^{2y}\/(4y)^{n\/2}}.
$$
 Then the last estimate gives 
\[
\lb{626}
|\sum_{n\ge 5}L_n(z)|\le e^{y}|\sum_{n\ge 5}\ve ^n|={\ve ^5e^{y}\/1-\ve}\le 2\ve ^5e^{y}.
\]
 \er{612}-\er{614} and  \er{617} give $L_s(z)=o(e^{2y})$
 as $y\to \iy$. Then combining  this and \er{626},\er{618}
 we obtain \er{619}.

Using the asymptotics $L=\cos z (I_N+S+O(|S|^2))$
and $\det (I+S)=\exp (\Tr S+O(|S|^2))$,
where $S={i\mV+o(1)\/2z}$ as $\Im z\to \iy$ we obtain \er{620}.
\BBox

\no {\bf Proof of Theorem \ref{T23}.} 
Recall the simple fact: Let $A, B$ be matrices and 
and $\s(B)$ be spectra of $B$. If $A$ be normal, then $\dist\{\s(A),\s(A+B)\}\le |B|$ (see p.291 \cite{Ka}).

 From \er{610} we have
 $L(z)=\cos zI_N+o(e^{|\Im z|}),\ |z|\to \iy $,
 where the operator $\cos zI_N$ has the eigenvalues $\cos z$ of the
multiplicity $N$. Due to the result from [Ka] and asymptotics
above we deduce that the eigenvalues $\D_m(z)$ of matrix $L(z)$
satisfy the asymptotics $\D_j(z)=\cos z+o(e^{|\Im z|}), j=1,..,N$
which gives that $M\in \mM_N^0$.

Lemma \ref{T62} implies $M\in \mM_N^{0,0}$.
Thus using  Theorem \ref{T1}-\ref{T2}, we obtain the proof of  Theorem \ref{T23}.
\BBox

\no {\bf 2. The Schr\"odinger operator.}
In order to prove Theorem \ref{T22} we determine the asymptotics of $M$. The fundamental solution $\vp$ satisfy the integral equations
\[
 \lb{627}
\vp(t,z)= \vp_0(t,z)+\int_0^t{\sin z(t-s)\/z}V(s)\vp(s,z)ds,\ \ \
\vp_0(t,z)={\sin z t\/z}I_N,\ \ \
\]
where $(t,z)\in\R\ts\C$. The standard iterations in \er{627} yields
\[
\lb{628}
 \vp(t,z)={\sum}_{n\ge 0} \vp_n(t,z)\,, \quad
\vp_{n+1}(t,z)= \int_0^t{\sin z(t-s)\/z}V(s)\vp_n(s,z)ds.
\]
The similar expansion $\vt={\sum}_{n\ge 0} \vt_n$ with
$\vt_0(t,z)=I_N\cos zt$ holds. We need

\begin{lemma} \lb{T63}
Let $V\in\mH$ and let $\vk={\|V\|\/|z|_1}$ and $|z|_1=\max\{1,|z|\}$. Then for any integers $m\ge t\ge 0, n_0\ge -1$ the following estimates hold:
\begin{multline}
\lb{629}
\max\lt\{\lt|\vt(1,z)-\sum_0^{n_0}\vt_n(1,z)\rt|,
\lt||z|_1\lt(\vp(1,z)-\sum_0^{n_0}\vp_n(1,z)\rt)\rt|,
\lt|{1\/|z|_1}\lt(\vt'(1,z)-\sum_0^{n_0}\vt_n'(1,z)\rt)\rt|,\\
\lt|\vp'(1,z)-\sum_0^{n_0}{\vp^{(n)}}'(1,z)\rt|\rt\} \le
{\vk^{n_0+1}\/(n_0+1)!}e^{(|\Im z|+\vk)},
\end{multline}
\[
\lb{630}
 |T_m(z)-2N\cos mz-{\sin mz\/z}m\Tr V^0|\le 2N m\vk e^{m(|\Im z|+\vk)}.
\]
\end{lemma}
\no {\bf Proof.} We prove the estimates of $\vp$, the proof for
$\vp',\vt,\vt'$ is similar. \er{628} gives
$$
\vp_n(t,z)= \int_{D_n}f_n(t,s)V(s_1)\cdot ..\cdot V(s_n)ds, \ \ \ \
\ f_n(t,s)=\vp_0(s_n,z)\prod_1^n {\sin z(s_{k-1}-s_{k})\/z},
$$
where $s=(s_1,..,s_n)\in \R^n, s_0=t$ and $D_n=\{0<s_n<
...<s_{2}<s_{1}<t\}$. 
Substituting the
estimate $|\vp_0(t,z)|=|z^{-1}\sin zt|\le |z|_1^{-1}e^{|\Im z|t}$ into the last integral, we obtain
$$
|\vp_n(t,z)|\le{e^{|\Im z|t}\/|z|_1^{n+1}}\int_{D_n}|V(s_1)|\cdot ..\cdot |V(s_n)|ds\le
{e^{|\Im z|t}\/|z|_1^{n+1}}\cdot{1\/n!}\lt(\int_0^t|V(x)|dx\rt)^n.
$$
This shows that for each $t\ge 0$
the series \er{628} converges uniformly on bounded subset of $\C$.
Each term of this series is an entire function. Hence the sum is
an entire function. Summing the majorants we obtain estimates
 \er{629}. 
 
 The function $T_m,m \ge 1$ is entire, since $\vp,\vt$
 are entire. We have 
 $T_m=\Tr M^m(z)=\sum_{n\ge 0}T_{m,n}(z)$, where
$$
 T_{m,0}=2N\cos mz,\ \ \ \ T_{m,n}(z)=\Tr
\vt_n(m,z)+\Tr \vp_n'(m,z),\ \ n\ge 1,
$$
$$
T_{m,1}(z)={1\/z}\int_0^m\!\!\! (\sin z(m-t)\cos zt+\cos
z(m-t)\sin zt )\Tr V(t)dt={\sin mz\/z}m\Tr V^0.
$$
The estimates $|\Tr \vp_n'(m,z)|\le {(m\vk)^n\/n!}e^{|\Im z|m}$ and
$|\Tr \vt_n(m,z)|\le {(m\vk)^n\/n!}e^{|\Im z|m}$ give
$
|T_{m,n}(z)|\le (2N){(m\vk)^n\/n!}e^{m|\Im z|},\ \ n\ge 0,
$
which yields \er{630}.
 $\BBox$

We will obtain the simple properties of the monodromy matrix.
We introduce the modified monodromy matrix $\P$ and 
the matrix $\L$ by
\[
\lb{631} \P=\mU^{-1} M \mU
=\ma \vt(1,z)& z\vp(1,z)\\ z^{-1}\vt'(1,z)& \vp'(1,z)\am,
 \ \ \ \mU=I_N\os zI_N \ \ \ z\in \C,
\]
\[
\lb{632}
\L=\mU^{-1} L\mU={\P+\P^{-1}\/2}={1\/2}\ma \vt(1,z)+\vp'(1,\ol z)^* & z(\vp(1,z)-\vp(1,\ol z)^*\\
z^{-1}(\vt'(1,z)-\vt'(1,\ol z)^*) & \vt(1,\ol z)^*+\vp'(1,z)\am,
\]
where we used the identity \er{cm-1}. Note that
$M$ and $\P$ have  the same eigenvalues and the same traces. 
Using \er{629} we obtain
\[
\lb{633} \L(z)=\L_1(z)+{\L_2(z)\/2z^2}+{\L_3(z)\/2z^3}+{O(e^{|\Im z|})\/z^{4}},\ \ \ \ \L_1(z)=\cos z+{\sin z\/2z}V^0 
\ \ as \ \ |z|\to \iy,
\]
where
\[
\lb{634} \L_2(z)=\int_0^1\!\! dt \int_0^t \sin z(t-s)
\ma a_{11}(t,s,z)& a_{12}(t,s,z)\\
a_{21}(t,s,z)& a_{22}(t,s,z)\am  ds, \ \ \ z\in \C,
\]
\[
\lb{635} a_{11}=\sin z(1-t)\cos zs V_tV_s+ \cos
z(1-t)\sin zs V_sV_t,\ \ \ a_{22}(z)=a_{11}(\ol z)^*,
\]
\[
\lb{636} a_{12}=\sin z(1-t)\sin zs (V_tV_s-V_sV_t), \ \
\  a_{21}=\cos z(1-t)\cos zs (V_tV_s-V_sV_t),
\]
where $V_t=V(t)$ and
\[
\lb{637} \L_3(z)=\int_0^1\!\! dt\!\!  \int_0^t\!\! ds\!\!
\int_0^s\!\!   \sin z(t-s)\sin z(s-u)
\ma b_{11}(t,s,u,z)& b_{12}(t,s,u,z)\\
b_{21}(t,s,u,z)& b_{22}(t,s,u,z)\am  du,
\]
\[
\lb{638} b_{11}=\sin z(1-t)\cos zu V_tV_sV_u+ \cos z(1-t)\sin
zu V_uV_sV_t,\ \ \ b_{22}(z)=b_{11}(\ol z)^*,
\]
\[
\lb{639} b_{12}=\sin z(1-t)\sin zu (V_tV_sV_u-V_uV_sV_t), \
\ \
b_{21}=\cos z(1-t)\cos zu (V_tV_sV_u-V_uV_sV_t).
\]

\begin{lemma} \lb{T64}
For each $(r,V)\in \R_+\ts\mH$ and $V^*=V$ the following asymptotics hold:
\[
\lb{640} 2^{2N}\det \L(z)=\exp \rt(-2Niz+i{\Tr V^0\/z}
 +{i\|V\|^2+o(1)\/4z^3}\rt),
\]
\[
\lb{641} \Tr \L_2(z)=\rt(-B_2+{i\|V\|^2+o(1)\/2z}\rt)\cos z, 
\]
\[
\lb{642} \Tr V^0\L_2(z)=-3B_3\cos z+O(e^{|\Im z|}/z),
\]
\[
\lb{643} \Tr \L_3(z)=-iB_3{\cos z\/2}+o(e^{|\Im z|})
\]
 as $y\ge r|x|, y\to \iy$, where  $ V^{0}=\!\!\int_0^1V(t)dt,\ \ B_{n}=\Tr{(V^{0})^n\/n!}$. 
\end{lemma}
\no {\it Proof.} Let $V_t=V(t)$.  Asymptotics \er{TA1}, \er{TA2} yield 
$
\int_0^1\int_0^t\!\cos z(1-2t+2s)\Tr V_tV_s={i\|V\|^2+o(1)\/2z}\cos z$. Then using \er{634} and $
\int_0^1\!\! \int_0^t\Tr V_tV_sdtds={1\/2}\Tr
\lt(\int_0^1V_tdt\rt)^2=B_2,
$ we obtain
$$
\Tr \L_2(z)=
2\int_0^1\int_0^t\!\!\sin z(t-s)\sin z(1-t+s)\Tr V_tV_sdtds
$$
$$
=\int_0^1\int_0^t\!\!\! (\cos z(1-2t+2s)-\cos
z)\Tr V_tV_sdtds={i\|V\|^2+o(1)\/2z}\cos z- B_2\cos z,
$$
which yields \er{641}.
We show \er{642}. Let $G(t,s)=\Tr V_0(V_tV_s+V_sV_t)$. Using
\er{634}, \er{TA1}, \er{TA2} we have
$$
\Tr V^0\L_2(z)=\int_0^1\int_0^t\!\!\! \sin z(t-s)z(1-t+s)
\Tr V^0(V_tV_s+V_sV_t)ds
$$
$$
={1\/2}\int^1dt\!\int_0^t\!\! (-\cos
z+\cos z(1-2t+2s))G(t,s)ds=-3B_3\cos z +O(e^{|\Im z|}/z),
$$
since
$$
\int_0^1dt\int_0^t\!\!\! \Tr V_0(V_tV_s+V_sV_t)ds
=\Tr V^0\lt(\int_0^1V_tdt\rt)^2=6B_3,
$$
which yields \er{642}.
Consider $\Tr \L_3$. The identity \er{637} gives
$$
\Tr \L_3(z)=\int_0^1dt\int_0^t\!\! ds\! \int_0^s\!\!\!
\sin z(t-s)\sin z(s-u)
\Tr (b_{11}+b_{22})du
$$
$$
=\int_0^1dt\int_0^t\!\!ds\! \int_0^s\!\!
\sin z(t-s)\sin z(s-u)\sin z(1-t+u)Rdu,\ \ \ \ \ \ R=\Tr (V_tV_sV_u+V_uV_sV_t).
$$
Using the identity
$$
4\sin z(t-s)\sin z(s-u)\sin z(1-t+u)=-\sin z+P,
$$$$P=\sin z(1-2s+2u)-\sin z(1-2t+2s)-\sin z(1-2t+2u),
$$
we get
$$
\L_3=-\L_3^0{\sin z\/2} +\L_3^1,\ \ \
\L_3^0={1\/2}\int_0^1dt\int_0^t\! ds\!\int_0^s\! Rdu,\ \
\L_3^1={1\/4}\int_0^1dt\int_0^t\! ds\!\int_0^s\! P Rdu,
$$
where
$$
\L_3^0={\Tr\/2}\int_0^1V_tdt\int_0^t\! ds\!\int_0^s\!
(V_sV_u+V_uV_s)du=
{\Tr\/2}\int_0^1V_t\rt(\int_0^t\! V_s ds\rt)^2
=B_3.
$$
Due to \er{TA-3} we obtain $\L_3^1=o(e^{|\Im z|})$, which yields \er{643}.

We will show \er{640}. Asymptotics \er{633} yields
\[
\lb{644}
{\L\/\cos z}=I_{2N}+S,\ \
S=i{V^0I_{2N}\/2z}+{\L_2\/2z^2\cos z}+{\L_3\/2z^3\cos z}+O(z^{-4}),\ \
\]
as $|z|\to \iy, y\ge r|x|$,
since $\sin z=i\cos z+O(e^{-y})$. In order to use the identity
$$
\det (I+S)=\exp \rt(\Tr S-\Tr {S^2\/2}+\Tr {S^3\/3}+o(z^{-3})\rt),\ \ |S|=O(1/z),
$$
we need the traces of $S^m, m=1,2,3$. Due to \er{641}-\er{643} we
get
\[
\lb{645}
{\Tr S^3\/3}=-i\Tr {(V^0)^3I_{2N}\/3(2z)^3}+O(z^{-4})
=-i {B_3\/2z^3}+O(z^{-4}),
\]
\[
\lb{646}
-\Tr{S^2\/2}={\Tr\/2}\rt({(V^0)^2\/(2z)^2}I_{2N}-
2i{V^0\L_2\/4z^3\cos z}+O(z^{-4})\rt)=
{B_2\/2z^2}+i{3B_3\/4z^3}+O(z^{-4})
\]
and
\[
\lb{647}
\Tr S=\Tr \rt(i{V^0\/2z}+{\L_2\/2z^2\cos z}+
{\L_3+O(z^{-1})\/2z^3\cos z}\rt)
=i{B_1\/z}+\rt(-{B_2\/2z^2}+i{\|V\|^2\/4z^3}\rt)-i{B_3+o(1)\/4z^3}
\]
and summing \er{645}-\er{647} we get \er{640}. 
$\BBox$

\no {\bf Proof of Theorem \ref{T22}.}  
Recall the simple fact: Let $A, B$ be matrices and 
and $\s(B)$ be spectra of $B$. If $A$ be normal, then $\dist\{\s(A),\s(A+B)\}\le |B|$ (see p.291 \cite{Ka}).

 From \er{629}, \er{633} we have
 $
 \L(z)=\L_1(z)+O(z^{-2}e^{|\Im z|}),\ |z|\to \iy
 $
 where the diagonal operator  $\L_1=\cos z+{\sin z\/2z}V^0$
has the eigenvalues $\D_j^0(z)=\cos z-V_{j}^0{\sin z \/2z}, j\in \ol{1, N}$ with the multiplicity 2. Using the result from [Ka] and asymptotics
above we deduce that the eigenvalues $\D_m(z)$ of matrix $\L(z)$
satisfy the asymptotics $\D_j(z)=\D_j^0(z)+O(z^{-2}e^{|\Im z|}),
 j\in \ol N$. Then $M\in\mM_N^0$ and  Lemma \ref{T64} yields 
$M\in\mM_N^{0,1}$.

The function $\F(iy,\n)$ is real for $y,\n\in \R$. Then $\f_j(z)=\ol \f_j(-\ol z)$ for all $(z,j)\in \ol\C_+\ts \ol {1, N}$. This yields that the set $\{\D_m(z)\}_1^N=\{\D_m(-\ol z)\}_1^N,z\in \ol\C_+$, which gives $q(z)=q(-\ol z), z\in \ol\C_+$.
Thus $q(x)=q(-x)$ for all $x\in \R$ and the identities
$$
k(-\ol z)=-\ol z+{1\/\pi }\int_\R{q(t)dt\/t+\ol z}=-\ol z-{1\/\pi }\int_\R{q(s)dt\/s-\ol z}=-\ol k(z), z\in \C_+
$$
 give $-k(-\ol z)=\ol k(z), z\in \C_+$.
 Thus by Theorem \ref{T1}-\ref{T2}, we obtain the proof of Theorem \ref{T22}  with the exception of \er{T22-6}, \er{T22-7}.
 Note that similar arguments give
$$
Q_4=I_2^D+I_4^S-Q_0Q_2=\int_0^1 {\Tr (V'(t)^2+2V^3(t))dt\/2^5N},
\ \ \ \ \ if \ \ V'\in \mH.
$$
 
 Defining $z_n^0={z_n^++z_n^-\/2}, r={|g_n|\/2}$
and using $(z_n^0+x)^{2}+(z_n^0-x)^{2}\ge 2(z_n^0)^{2}$, we have 
$$
\int_{g_n}{t^{2}q(t)dt\/\pi}\ge
\int_{g_n}{t^{2}q_n^0(t)dt\/\pi}\ge \int_0^r((z_n^0+x)^{2}+(z_n^0-x)^{2}){\sqrt
{r^2-x^2}\/\pi}dx\ge
{r^2\/2}(z_n^0)^{2}={|\g_n|^2\/32},
$$
$$
Q_2={1\/\pi}\int_\R t^2q(t)dt=\sum_{n\in \Z}{1\/\pi}\int_{g_n} t^2q(t)dt
\ge {1\/32} \sum_{n\in \Z}
|\g_n|^2={1\/16} \sum_{n\ge 1} |\g_n|^2,
$$
 which yields  \er{T3-5}.
 
 Assume that $\s_{(N)}=\s(F)$. In this case the function
 $h$ is given by: 
\[
\lb{deh}
 h(p)=0, p\ne {\pi\/N} \Z,\  \ \ \ \ {\rm and}\ \ \ \ 
 (h({\pi\/N}n))_{n\in \Z}\in \ell_1^2.
\] 
Recall estimates from \cite{K9}.
Consider the conformal mapping $k:\C_+\to \K(h)$ for the case $k(z)=z+o(1)$ as $|z|\to \iy$ and $h$ given by \er{deh} \cite{K9}.
Then for some absolute constant $C_0$ the following estimate was obtained: $Q_0\le C_0G_0^2(1+G_0^2)$ and 
$Q_2\le C_0G_2^2(1+G_2^{1\/3})$, which  gives the second estimates  \er{T3-4} and \er{T3-5}.
\BBox

\section {Appendix}
\setcounter{equation}{0}

\begin{lemma} \lb{TA1}
Let functions $h_1,..h_n\in L^2(0,1)$ for some $n\ge 3$. 
Then
\[
\lb{TA-1} \int_0^1dt\int_0^t e^{i2z(t-s)}h_1(t)h_2(s)ds ={i\/2z}\lt(\int_0^1h_1(t)h_2(t)dt+o(1)\rt)
\]
\[
\lb{TA-2} \int_0^1dt\int_0^t
e^{-i2z(t-s)}h_1(t)h_2(s)ds={o(e^{2y})\/|z|},
\]
\[
\lb{TA-3} \int_0^1dt\int_0^tds\int_0^s e^{\pm iz(1-2\z)}h_1(t)h_2(s)h_3(u)du=o(e^{y}),
\]
\[
\lb{TA-4} 
\int_0^1dt_1...\int_0^{t_{n-1}}dt_n
e^{\pm iz(1-2t_1+2t_2\dots +(-1)^n2t_n)}\prod_1^n h_j(t_j)
={o(e^{y})\/|z|},\ \ all \ \ n\ge 3,
\]
as $r|x|<y\to \iy$ for any fixed $r>0$, where $\z$ is one of functions: $s-u, t-u$, or $t-s$.
\end{lemma}

\no {\it Proof.} Let $F(t,s)=h_1(t)h_2(s), t,s\in (0,1)$. We have 
\[
\lb{TA-5} f_1(z)\ev\int_0^1dt\int_0^te^{i2z(t-s)}F(t,s)ds
={1\/2}\int_0^1ds\int_0^1e^{i2z|t-s|}F(t,s)ds.
\]
  Substituting the identities
$$
{2z\/ \pi i}\int_\R  {e^{i(t-s)k}dk\/k^2-4z^2}=e^{i2z|t-s|},\ \
\ z\in \C_+,\ \ \ \ 
\hat F(k)\ev{1\/ 2\pi }\iint_{[0,1]^2}
e^{ik(t-s)}F(t,s)dtds, 
$$
into \er{TA-5} we obtain
$$
f_1={2z\/i}\int_\R {\hat F(k)dk\/k^2-4z^2}
={1\/i2z}\int_\R \rt(-1+{k^2\/k^2-4z^2}\rt)\hat F(k)dk=
{i\/2z}\rt(\int_0^1 F(t,t)dt+o(1)\rt).
$$
which yields \er{TA1}. Consider $f(z)\ev\int_0^1dt\int_0^te^{i2z(t-s)}F(t,s)ds$. 
We have
\[
\lb{TA-6} |f(z)|^2\le g(z)\|h_1\|^2\|h_2\|^2,\ \  \ \ g(z)=\int_0^1dt\int_0^t
e^{4y(t-s)}ds\le \int_0^1 {e^{4yt}\/4y}dt\le {e^{4y}\/(4y)^2}.
\]
Let $F_0$ be a smooth function such that $\|F-F_0\|=\ve$ for
some small $\ve>0$. Define the function
$f_0(z)=\int_0^1dt\int_0^t e^{-i2z(t-s)}F_0(t,s)ds$. Using
\er{TA-6} we obtain
\[
\lb{TA-7}
|f(z)|\le |f_0(z)|+|f(z)-f_0(z)|\le |f_0(z)|+
\|F-F_0\|{e^{2y}\/(4y)}
\]
and the integration by parts yields $f_0(z)=O({e^{2y}\/y^2})$.
Thus we obtain \er{TA-2}, since $\ve $ is arbitrary small.
 The similar arguments yield \er{TA-3} and \er{TA-4}.
$\BBox$

We formulate the following results from [CK].

\begin{lemma} \lb{TA2}
The function $f(z)=\log |\x(z)|, z\in \C\sm [-1,1]$ is subharmonic
and continuous in $\C$. Moreover, for some absolute constant $C$
the following estimate is fulfilled:
\[
c |f(z)-f(z_0)|\le C\ve^{1\/2},\ \ \ if \ \ |z-z_0|\le \ve
\max\{2,|z_0|\}, \ \ 0\le\ve\le {1\/8}, \ \ z,z_0\in \C.
\]
\end{lemma}

 \no {\bf Acknowledgments.}
The author was partly supported by DFG project BR691/23-1.
The various parts of this paper were written at the Mittag-Leffler Institute, Stockholm  and in the Erwin Schr\"odinger Institute for Mathematical Physics, Vienna, the author is grateful to the Institutes for the hospitality.


\end{document}